\documentclass[11pt]{article}
\usepackage[dvips]{graphicx}
\usepackage{cite}
\usepackage{amsmath,amssymb,subfigure,epsfig,bbm,stmaryrd,MnSymbol,mathrsfs,url}
\usepackage{pstricks}
\usepackage{subfigure}
\usepackage{booktabs}
\usepackage{multirow}
\usepackage{amsmath,amssymb}
\usepackage{pstricks}
\usepackage{amsmath,amssymb,subfigure,epsfig,bbm,stmaryrd,MnSymbol,mathrsfs}
\usepackage{subfigure}
\DeclareMathAlphabet{\mathpzc}{OT1}{pzc}{m}{it}
\DeclareFontFamily{OT1}{pzc}{}
\DeclareFontShape{OT1}{pzc}{m}{it}{<-> s * [1.2] pzcmi7t}{}
\DeclareMathAlphabet{\mathpzc}{OT1}{pzc}{m}{it}

\usepackage{algorithm,mdframed}
\usepackage{algorithmic}

\newcommand{\bdelta}{\boldsymbol{\delta}}

\newcommand{\balpha}{\boldsymbol{\alpha}}
\newcommand{\je}{\boldsymbol{\mathcal{J}}_{\!\mathcal{E}}}
\newcommand{\he}{\boldsymbol{\mathcal{H}}_{\mathcal{E}}}
\newcommand{\ii}{\mathfrak{i}\;}

\begin{document}
\newcommand{\BigFig}[1]{\parbox{12pt}{\Huge #1}}
\newcommand{\BigZero}{\BigFig{0}}
\title{Sparse Shape Reconstruction}   

\author{Alireza Aghasi and Justin Romberg\thanks{School of Electrical and Computer Engineering, Georgia Institute of Technology, Atlanta, GA. Emails: {\tt aaghasi@ece.gatech.edu} and {\tt jrom@ece.gatech.edu}.}}

\date{}    
\maketitle

\begin{abstract}
This paper introduces a new shape-based image reconstruction technique applicable to a large class of imaging problems formulated in a variational sense. Given a collection of shape priors (a shape dictionary), we define our problem as choosing the right elements and geometrically composing them through basic set operations to characterize desired regions in the image. This combinatorial problem can be relaxed and then solved using classical descent methods. The main component of this relaxation is forming certain compactly supported functions which we call ``knolls'', and reformulating the shape representation as a basis expansion in terms of such functions. To select suitable elements of the dictionary, our problem ultimately reduces to solving a nonlinear program with sparsity constraints. We provide a new sparse nonlinear reconstruction technique to approach this problem. The performance of proposed technique is demonstrated with some standard imaging problems including image segmentation, X-ray tomography and diffusive tomography.
\end{abstract}


\section{Introduction}
In many imaging applications, the main objective is to identify and characterize regions of interest in a given domain. This characterization is usually based on some property defined over the imaging domain. For instance in image segmentation \cite{pham2000current, pal1993review}, this property is directly or statistically related to the pixel values and the partitioning is usually meant to aggregate similar regions. In shape-based inverse problems \cite{kirsch1998characterization, chung2005electrical, dorn2000shape}, the property of interest corresponds to a spatial physical parameter which needs to be determined by some indirect observations. Here a shape-based characterization delineates inclusions and obstacles causing a contrast in the values of the spatial parameter. In the basic binary case the shape characterization problem is formulated as partitioning a compact imaging domain $D\in \mathbb{R}^q$, into $\tilde \Omega$ and $D\setminus \tilde \Omega$, where $\tilde \Omega\in \mathbb{R}^q$ is a closed set representing the object(s) of interest.

In this paper we propose a new approach to shape-based modeling. The idea is to characterize the target inclusion (or partition) as a composition of given shape prototypes. These prototypes may either be shapes with simple geometries or shape priors that are likely to be related to the structure of the target geometry. The idea of composing basic shapes to form more complex structures is quite intuitive. The challenge, is to formulate shape composition concretely and then incorporate it into the reconstruction process in a computationally tractable manner.

For the purpose of shape composition, we consider a mechanism to merge shapes and exclude undesirable regions from the aggregate to generate $\tilde \Omega$. More formally, given a set of prototype shapes, $\mathcal{S}_1, \mathcal{S}_2, \cdots \mathcal{S}_{n}$, which are known closed sets in $\mathbb{R}^q$,
we formulate the shape composition process as applying set union and relative complement among a selected number of shapes in the reference set to generate a region that approximates $\tilde \Omega$. A brute force approach would entail exploring all possible selections and composition possibilities to find a good fit. To have an idea about the accuracy of an estimate of $\tilde \Omega$, we specifically focus on variational shape-based problems where a score (value of an energy functional) is assigned to every closed region $\Omega\in D$ and $\tilde \Omega$ is meant to minimize such functional.

This type of brute force search is in general computationally intractable. In this paper we will propose a relaxation formed by defining so called \emph{knoll} functions which characterize $\mathcal{S}_1, \mathcal{S}_2, \cdots \mathcal{S}_{n}$ by their support. We will show that for a smooth energy functional associated with the shape-based problem, this relaxation can convert the problem into minimizing a smooth nonlinear function in $\mathbb{R}^n$.

To encourage our target shape to be as simple as possible, we impose sparsity constraints on the resulting nonlinear cost. Accordingly, we propose a minimization technique inspired by the idea developed by van den Berg and Friedlander in \cite{van2008probing}. The original idea in \cite{van2008probing} developed for linear inverse problems is merely applicable to the corresponding quadratic least squares cost. What we will put forth is a generalization of this approach applied to quadratic estimates of the nonlinear problem at iterative stages of approaching a minima.

The organization of this paper is as follows. In the remainder of the introduction section we provide a rather extensive overview of variational shape-based methods, specially in context of active contours, pixel based and parametric level set techniques. After providing this background, in Section 2 we present the shape composition idea and discuss employing a so called \emph{pseudo-logical} shape interaction property which will assist us in developing a relaxation to the corresponding combinatorial problem. Section 3 is devoted to developing a Gauss-Newton type sparsity promoting algorithm to solve the resulting nonlinear problem. Finally, in Section 4 we consider applying the proposed technique to some image segmentation and shape-based inverse problems\footnote{The Matlab code for the proposed algorithm is available at \url{http://users.ece.gatech.edu/aaghasi3/software.html}}. The image processing applications considered are segmentation with missing pixels and machine text recognition. We also consider an example of medical X-ray computed tomography and an archaeologic resistance tomography problem, both with limited available data, where majority of the state of the art techniques would not be able to provide satisfying reconstructions.

\subsection{Background on Variational Shape Reconstruction}
In the past decades considerable effort has been devoted to exploring variational techniques applicable to the shape-based characterization problem \cite{chan2003variational, kass1988snakes, chan2001active, osher2001level, dorn2006level, tsai2001curve}. The main intuition behind a variational approach is forming an energy functional, minimizing which tends to solve the characterization problem. More specifically, for an arbitrary bounded set $\Omega\in D$,
\begin{equation}\label{eq2}
\tilde \Omega=\operatorname*{arg\,min}_{\Omega} \mathpzc{E}( \Omega),
\end{equation}
where $\mathpzc{E}$ is the underlying energy or cost functional. An example of $\mathpzc{E}$ in the context of inverse problems is the data-model mismatch functional
\begin{equation}\label{eq3}
\mathpzc{E}(\Omega)=\|v-\mathcal{M}(\Omega)\|_{\mathbb{S}_v},
\end{equation}
where $v$ represents the observed data, $\mathcal{M}$ denotes the model mapping the shape geometry to the observations and $\mathbb{S}_v$ is a Hilbert space associated with the data. One of the main motivations in representing problems in variational forms is the chance of using descent optimization techniques.

Of course a straight-forward approach in determining a region $\Omega$ is determining its boundary $\mathcal{C}$. An early technique of this type is the active contour model (snakes), where the shape determination amounts to determining the parameters associated with a spline representation of $\mathcal{C}$ \cite{kass1988snakes}. Starting with an initial shape representation, this process takes steps along the descent direction of the energy to evolve the shape towards an optimum state. More specifically, by defining an artificial time $t$, an initial contour $\mathcal{C} (0)$ is evolved in time and according to
\begin{equation} \label{eq6}
\frac{\partial \mathcal{C}(t)}{\partial t}=-\;\mathpzc{E}'(\mathcal{C})\big[ \mathcal{C}(t)\big],
\end{equation}
to find the boundary of the shape that locally minimizes $\mathpzc{E}$ \cite{caselles1997geodesic}. Here $\mathpzc{E}'(\mathcal{C})[.]$ is a linear operator representing the G\^{a}teaux derivative (or first order variations) of the energy functional with respect to $\mathcal{C}$.

Another well known shape-based technique is the level set method \cite{osher1988fronts, osher2003level}. The main advantage of level sets over earlier methods such as snakes is their topological flexibility, dispelling the need to any prior assumptions about the number of connected components in $\tilde \Omega$. Here the zero level set of a Lipschitz continuous function, $\phi$, is used to identify $\mathcal{C}$. More specifically the objective is to determine $\phi$ such that
\begin{equation}\label{eq4}
  \left\{
     \begin{array}{ll}
       \phi(x)\geq 0 &  x \in \tilde \Omega \\
       \phi(x)<0 &  x \in D \setminus \tilde \Omega
     \end{array}.
  \right.
\end{equation}
By using a map as (\ref{eq4}) our geometric problem is cast as the calculus of variations problem
\begin{equation}\label{eq5}
\tilde \phi=\operatorname*{arg\,min}_{\phi \;\in \;\Theta} \mathpzc{E}(\Omega_{\phi}),
\end{equation}
where $\Omega_{\phi}$ is the shape resulted at the zero level set of $\phi$, and $\Theta$ represents a certain function space to prevent ill-conditioning. Most level set implementations consider elements of $\Theta$ to be signed distance functions (SDFs) \cite{osher2003level}.

For the level set methods, minimization of (\ref{eq5}) is performed by evolving an initial level set function, $\phi(x,t)=\phi_0(x)$, through the Hamilton-Jacobi equation
\begin{equation}\label{eq7}
\frac{\partial \phi(x,t)}{\partial t}+V(x,t)\cdot \nabla \phi(x,t)=0.
\end{equation}
This equation results from a straightforward differentiation of the front equation $\phi(x,t)=0$ with respect to $t$. In this equation $V(x,t)=\mbox{d}x/\mbox{d}t$ is a speed function applied to the zero level set of $\phi(x,t)$ and taken in the descent direction of $\mathpzc{E}$ to reduce it as $\phi$ evolves. At every iteration, through a proper speed function extension or re-initialization, the level set function is assured to remain an SDF. We refer the reader to \cite{osher2003level, dorn2006level, burger2005survey} for more details about implementation of this technique.

Level set methods are among the most successful techniques in shape representations, mainly due to their topological flexibility. They perform remarkable for a large class of image processing applications \cite{osher2003geometric, cremers2007review}, however, in general there are implementation concerns and performance deficiencies associated with them. Numerically speaking, evolving the level set function requires discretizing it over a dense grid of pixels and updating the pixel values over discrete time frames (iterations) \cite{osher2003level}. Using pixels to parameterize the level set function brings a large dimensionality to the problem which may severely affect the overall performance of the method in dealing with ill-posed functionals and inverse problems (see examples in \cite{aghasi2011parametric}). Although in principle level sets are topologically flexible, it is usually not possible to create holes or shapes away from the boundaries of the evolving shape \cite{burger2004incorporating}. Moreover, implementation complexities such as speed function extension and maintaining $\phi$ as an SDF are usually the inevitable components of this technique.

To overcome these problems and still take advantage of the topological flexibility, researchers have started considering parametric forms for the level set function \cite{aghasi2011parametric, bernard2009variational}. This strategy not only reduces problem's dimensionality by making it parametric, also provides a predetermined functional form to dispel the need to re-initialization. An example of such effort is the work by Bernard \emph{et al.} \cite{bernard2009variational}, where according to (\ref{eq5}), $\Theta$ is the space spanned by pre-assigned B-spline functions. The problem in this case usually amounts to the energy minimization
\begin{equation} \label{eq8}
\tilde{ \boldsymbol{\alpha}}=\operatorname*{arg\,min}_{\boldsymbol{\alpha}} \mathpzc{E}(\Omega_{ \phi(x,\boldsymbol{ \alpha})}),
\end{equation}
where $\boldsymbol{\alpha}$ is the vector of parameters associated with the parametric level set function $\phi(x,\boldsymbol{\alpha})$.

A main advantage of parameterizing the problem as (\ref{eq8}) is casting the shape characterization problem as a classic finite dimensional minimization problem, where gradient descent or even Newton type methods may be applied. For sufficiently smooth functionals, gradient and Hessian of $\mathpzc{E}$ with respect to elements of $\boldsymbol{\alpha}$ may be conveniently calculated by following the chain rule in Fr\'{e}chet spaces \cite{aghasi2011parametric}. More specifically
\begin{equation}\label{eq8a}
 \frac{\partial \mathpzc{E}}{\partial \alpha_i}=\mathpzc{E}'(\phi)[\frac{\partial \phi}{\partial \alpha_i}],
\end{equation}
and
\begin{equation}\label{eq8b}
 \frac{\partial^2 \mathpzc{E}}{\partial \alpha_i\partial \alpha_j}=\mathpzc{E}''(\phi)[\frac{\partial \phi}{\partial \alpha_i},\frac{\partial \phi}{\partial \alpha_j}]+\mathpzc{E}'(\phi)[\frac{\partial^2 \phi}{\partial \alpha_i\partial \alpha_j}],
\end{equation}
where the linear and bilinear operators $\mathpzc{E}'(\phi)[.]$ and $\mathpzc{E}''(\phi)[.,.]$ are respectively the first and second G\^{a}teaux derivatives of $\mathpzc{E}$ with respect to $\phi$. We refer the reader to \cite{aghasi2011parametric} for more details on a parametric level set representation.

Although a parametric form solves many concerns with traditional level sets, the problem still remains on choosing a suitable parametric form (e.g., a suitable set of basis functions) and efficiently determining the number of terms. Clearly a dense basis set to increase the resolution of the reconstructions may bring redundancy and ill-conditioning to the problem as the case with pixel based level sets.

In the sequel, we propose a new way of approaching the shape problem. We extend one of the basic ideas developed by Aghasi \emph{et al.} in \cite{aghasi2011parametric} and provide a new shape reconstruction technique applicable to a large class of problems. The method ultimately finds a parametric level set form, however, we link it to a systematic algorithm of choosing the right parameters in the course of reconstruction.

\section{Basic Shape Composition Idea}

In solving (\ref{eq2}) for the optimum shape $\Omega$, an intuitive idea would be to consider a collection of fixed shapes and reconstruct $\Omega$ by applying basic set operations on them. More specifically, consider the shape collection $\mathfrak{D}=\{\mathcal{S}_1, \mathcal{S}_2, \cdots \mathcal{S}_{n_d}\}$, where for $i=1,2, \cdots , n_d$, elements $\mathcal{S}_i$ represent closed known regions (shapes) in $\mathcal{D}$. We name $\mathfrak{D}$ as the \emph{shape dictionary} or simply the dictionary. Suppose through an oracle we know what the true shape $\Omega$ is. To logically (in the sense of set operations) express $\Omega$ in terms of elements of $\mathfrak{D}$, one way of approaching the problem is to select an appropriate set of shapes indexed by $\mathcal{I}_\oplus \subseteq \{1,2,\cdots, n_d\}$ and apply the union operation over them to form the more bulky superset $\bigcup_{i\in \mathcal{I}_\oplus} \mathcal{S}_i$ for $\Omega$. We now start to choose shapes indexed by $\mathcal{I}_\ominus \subset \{1,2,\cdots, n_d\}$ to carve out portions of $\bigcup_{i\in \mathcal{I}_\oplus}\mathcal{S}_i$ and make it a better approximation to $\Omega$ (in the sense of reducing $\mathpzc{E}$). To be concise, for a given dictionary $\mathfrak{D}$ we define the objective as the combinatorial problem of searching among elements $\Omega$ of the form
\begin{equation}\label{eq9}
 \Omega_{\mathcal{I}_\oplus,\mathcal{I}_\ominus}\triangleq \Big(\bigcup_{i\in \mathcal{I}_\oplus} \mathcal{S}_i \Big ) \setminus \Big(\bigcup_{j\in \mathcal{I}_\ominus } \mathcal{S}_j\Big),
\end{equation}
and find the suitable index sets $\mathcal{I}_\oplus $ and $\mathcal{I}_\ominus$ that minimize $\mathpzc{E}(\Omega)$.

This simple idea is inspired by approximation theory of suitably expressing a function as a linear combination of some given basis functions. Here instead of adding and subtracting the basis terms with suitable weights, suitable elements of $\mathfrak{D}$ are combined through the union and relative complement to provide the shape approximation.

We would like to highlight that more complex reference forms other than (\ref{eq9}) may be considered. We however maintain simplicity by suggesting this form and assuming that $\mathfrak{D}$ is rich enough that (\ref{eq9}) still provides the desired flexibility in shape representation. For instance if the dictionary is poor, consisting of only two shapes $\mathcal{S}_1$ and $\mathcal{S}_2$, and the true shape is $\tilde \Omega=\mathcal{S}_1 \cap \mathcal{S}_2$, none of the possibilities in the form of (\ref{eq9}) would be able to express $\tilde \Omega$. However if $\tilde \Omega$ is already among the elements of the dictionary or we consider a richer dictionary consisting of $\mathcal{S}_1$, $\mathcal{S}_2$, $\mathcal{S}_3=\mathcal{S}_1\setminus \mathcal{S}_2$ and $\mathcal{S}_4=\mathcal{S}_2\setminus \mathcal{S}_1$, we can find $\tilde \Omega$ by exploring different possibilities of the form (\ref{eq9}) since $\Omega=(\mathcal{S}_1\cup \mathcal{S}_2)\setminus (\mathcal{S}_3\cup \mathcal{S}_4)$.

\begin{figure}[!htbp]%
\hspace{-5mm}
\subfigure[][]{\includegraphics[width=67mm]{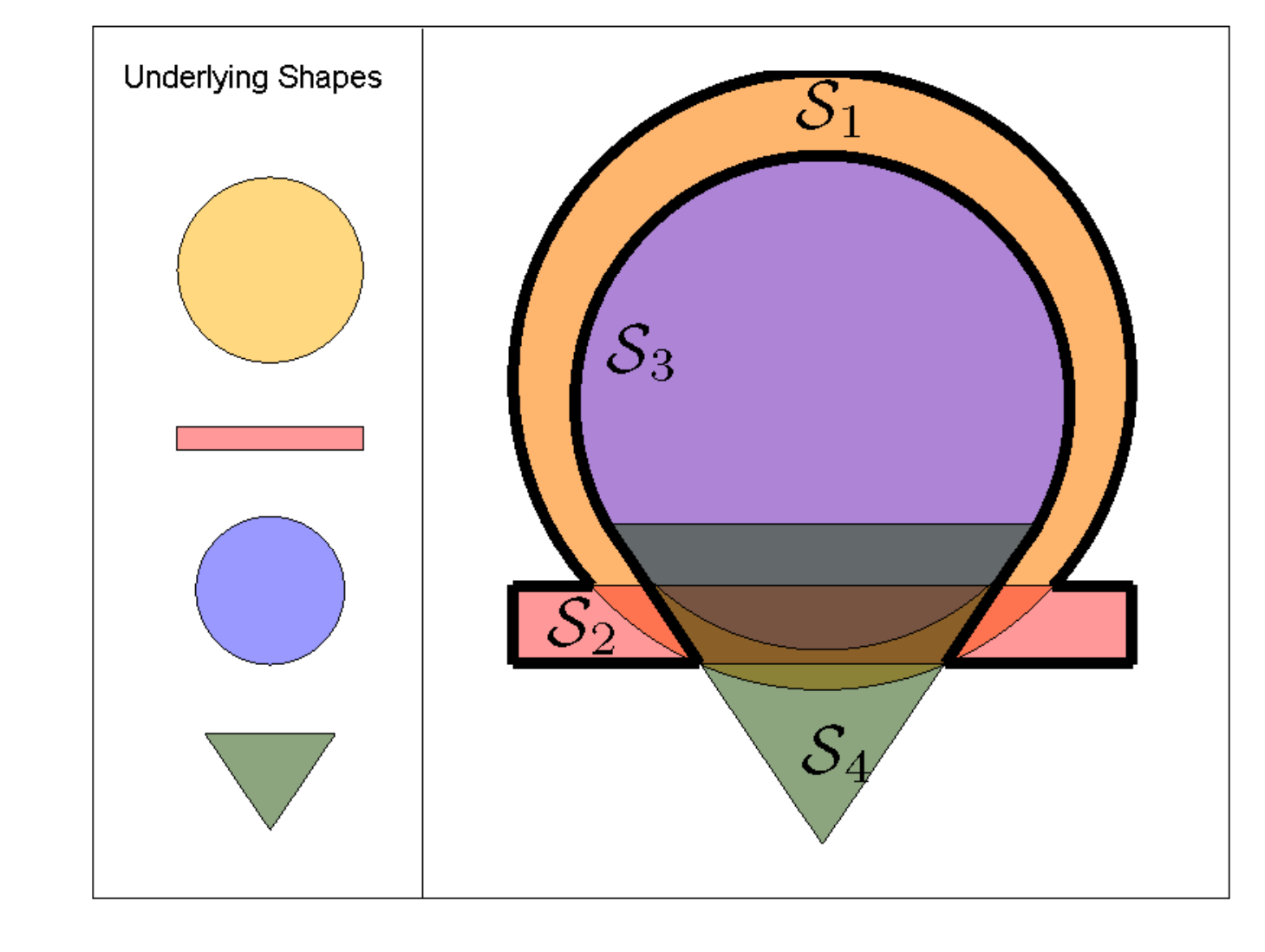}}
\subfigure[][]{\includegraphics[width=68mm]{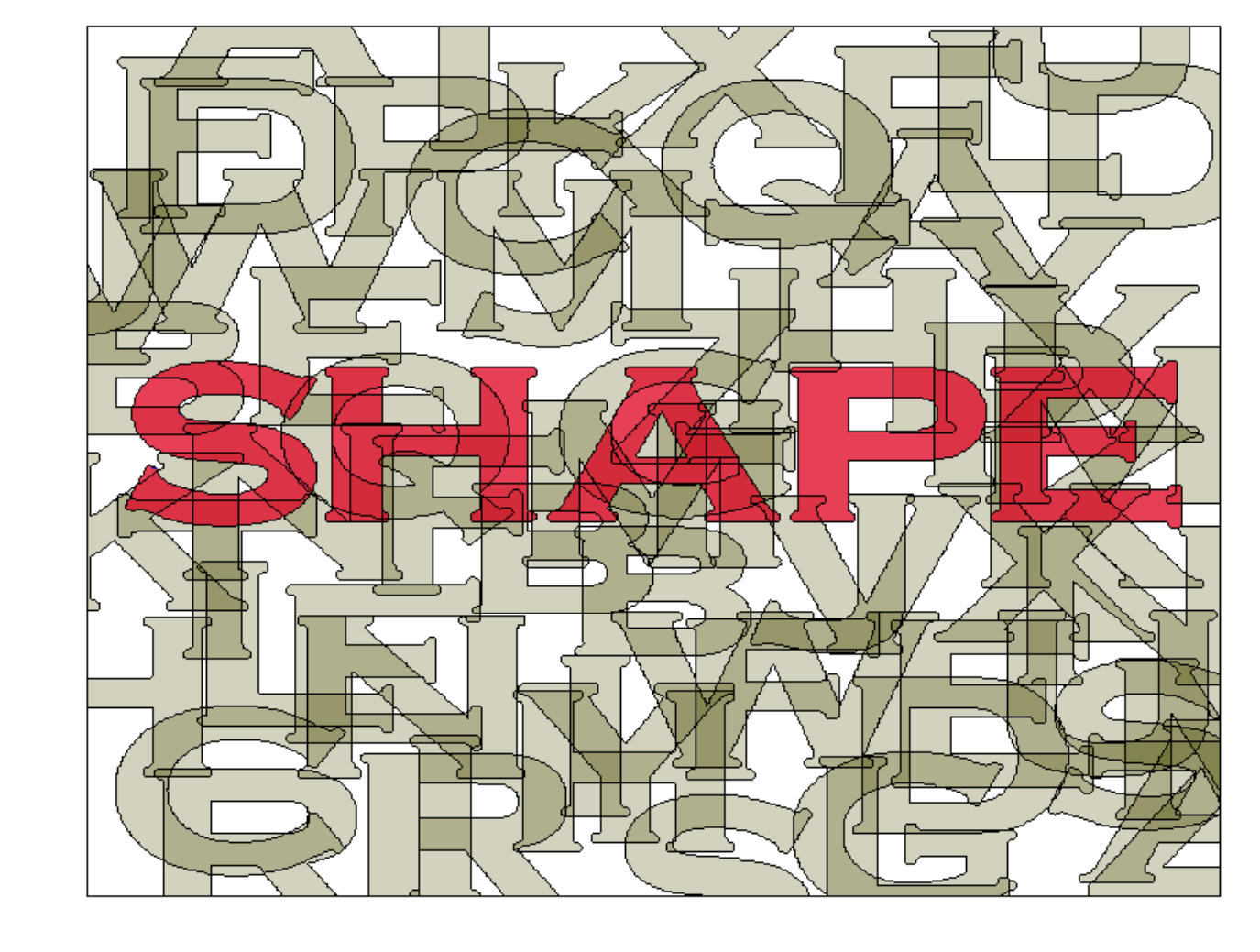}}
\caption{(a) Reconstruction of an ``Omega''-shaped region by applying basic set operations on simpler geometries (two circles, a rectangle and a triangle): in this case the desired shape can be written as $\Omega=(\mathcal{S}_1\cup\mathcal{S}_2)\setminus (\mathcal{S}_3\cup\mathcal{S}_4)$ (b) Choosing the right components of the desired shape among a dictionary of candidates: in this case extracting the word ``SHAPE'' from a collection of randomly placed characters}%
\label{fig1}%
\end{figure}

Selecting appropriate elements for $\mathfrak{D}$ may be based on the level of prior information about the geometric features of $\Omega$. The shapes $\mathcal{S}_i$ may be simple geometries that combine to form a more complex structure (e.g. see Fig \ref{fig1}(a)); a large collection of shape possibilities among which the true shape needs to be determined (Fig \ref{fig1}(b)); or a combination of both cases.

To bring this idea into application we need to develop a computationally tractable algorithm that performs the search among possible shape combinations and minimizes the functional in (\ref{eq2}) using a sufficiently sparse set of elements in the dictionary. More formally we reformulate (\ref{eq2}) as
\begin{equation}\label{eq10}
\{\mathcal{I}_\oplus, \mathcal{I}_\ominus\}=\operatorname*{arg\,min}_{|\mathcal{I}_\oplus|+|\mathcal{I}_\ominus|\leq s} \mathpzc{E}( \Omega_{\mathcal{I}_\oplus,\mathcal{I}_\ominus}),
\end{equation}
where $s$ represents the desired level of shape sparsity in the reconstructions. To provide an approximate solution to the combinatorial minimization problem (\ref{eq10}), in the sequel we present a relaxation strategy to maintain sufficient smoothness of the cost for gradient type optimization techniques to become applicable. We then provide a sparse reconstruction algorithm associated with nonlinear costs to employ a sparse number of dictionary elements in the final shape representation.

\subsection{Pseudo-Logical Shape Interaction}
In a recent work, Aghasi \emph{et al.} provided a rather general view of parametric level set methods for inverse problems. As an example of a parametric form, they proposed a basis expansion of smooth compactly supported radial basis functions (named as bumps) with adaptive dilation and center points. They brought into attention a so called ``pseudo-logical'' property of this class of functions that we generalize its notion to broader applicability and something not specific to smooth radial bumps.

Consider a given shape $\mathcal{S} \subset D\subset\mathbb{R}^q$. We use the terminology ``\emph{knoll}'' for a Lipschitz continuous function $\psi_\mathcal{S}: D\rightarrow [0,\infty)$ such that
\begin{equation}\label{eq11}
  \left\{
     \begin{array}{ll}
       \psi_\mathcal{S}(x)>0 &  x \in \mbox{int}(\mathcal{S}) \\
       \psi_\mathcal{S}(x)=0 &  x \in D\setminus \mathcal{S}
     \end{array},
  \right.
\end{equation}
with $\mbox{int}(\mathcal{S})$ denoting the interior of $\mathcal{S}$. Intuitively, a knoll $\psi_\mathcal{S}$ takes positive values inside $\mathcal{S}$ and vanishes outside $\mathcal{S}$ (and on the boundaries as a result of continuity). It can be easily inferred that for two knolls $\psi_{\mathcal{S}_1}$ and $\psi_{\mathcal{S}_2}$ corresponding to shapes $\mathcal{S}_1$ and $\mathcal{S}_2$, the following facts are true:
\begin{equation}\label{eq12}
\mbox{supp}(\psi_{\mathcal{S}_1}+\psi_{\mathcal{S}_2})=\mathcal{S}_1\cup \mathcal{S}_2
\end{equation}
and
\begin{equation}\label{eq13}
\lim_{\alpha \to +\infty} \mathscr{D}\big(\mbox{supp}^+(\psi_{\mathcal{S}_1}-\alpha \psi_{\mathcal{S}_2}), \mathcal{S}_1\setminus \mathcal{S}_2\big)=0.
\end{equation}
Here $\mbox{supp}^+(.)$ denotes the positive support, where the function takes values greater than zero and $\mathscr{D}(.,.)$ is a measure of dissimilarity between two shapes which basically vanishes for identical shapes (see \cite{cremers2006kernel} for examples of this measure). The basic message is summing up two knolls would imply a union operation on their supports, and subtracting a knoll of large weight $\alpha\gg 1$ from another knoll approximates applying relative complement on their positive supports.

\begin{figure}[!htbp]
\hspace{-.3cm}
\subfigure[][]{\includegraphics[width=65mm]{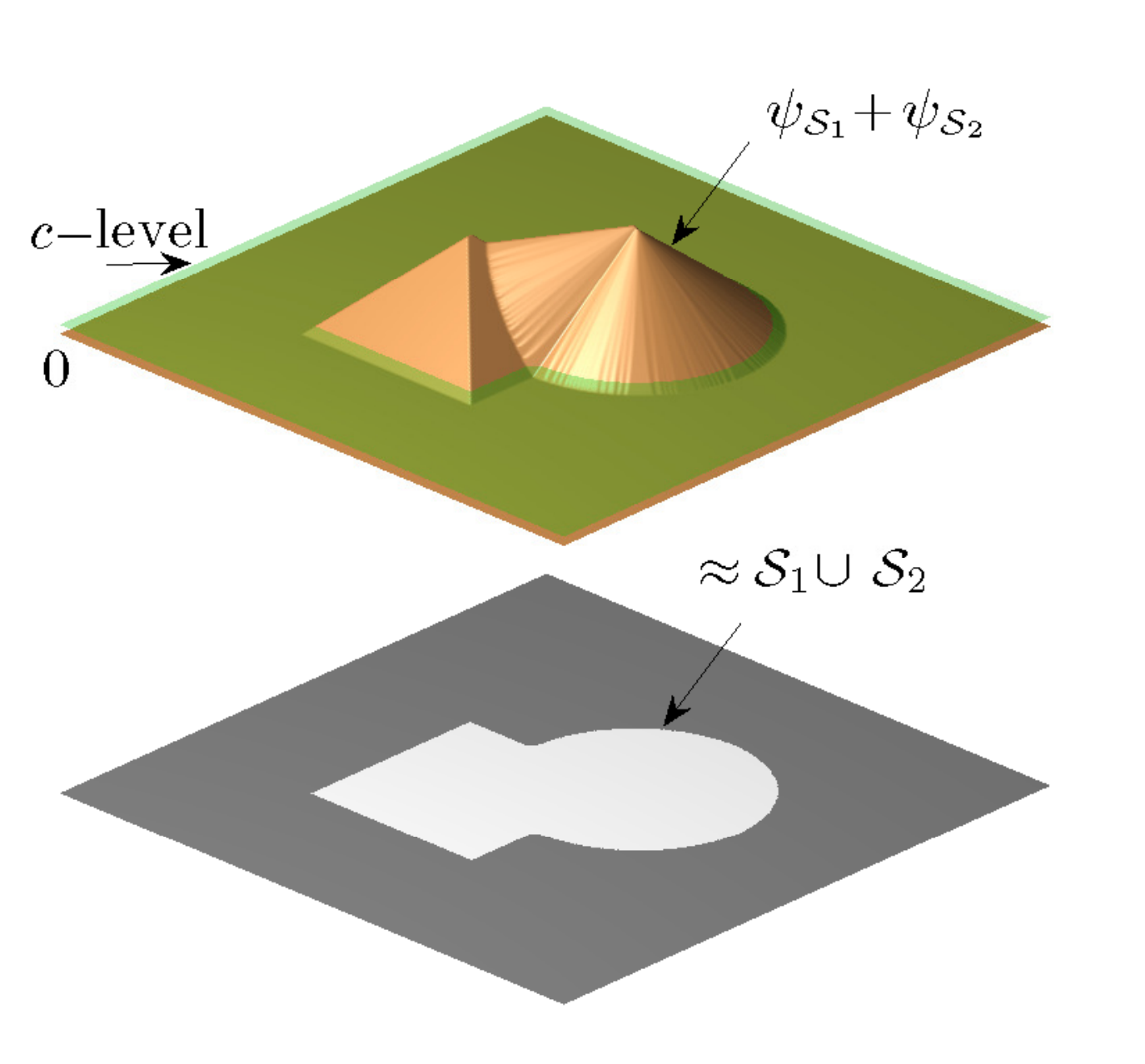}}
\subfigure[][]{\includegraphics[width=65mm]{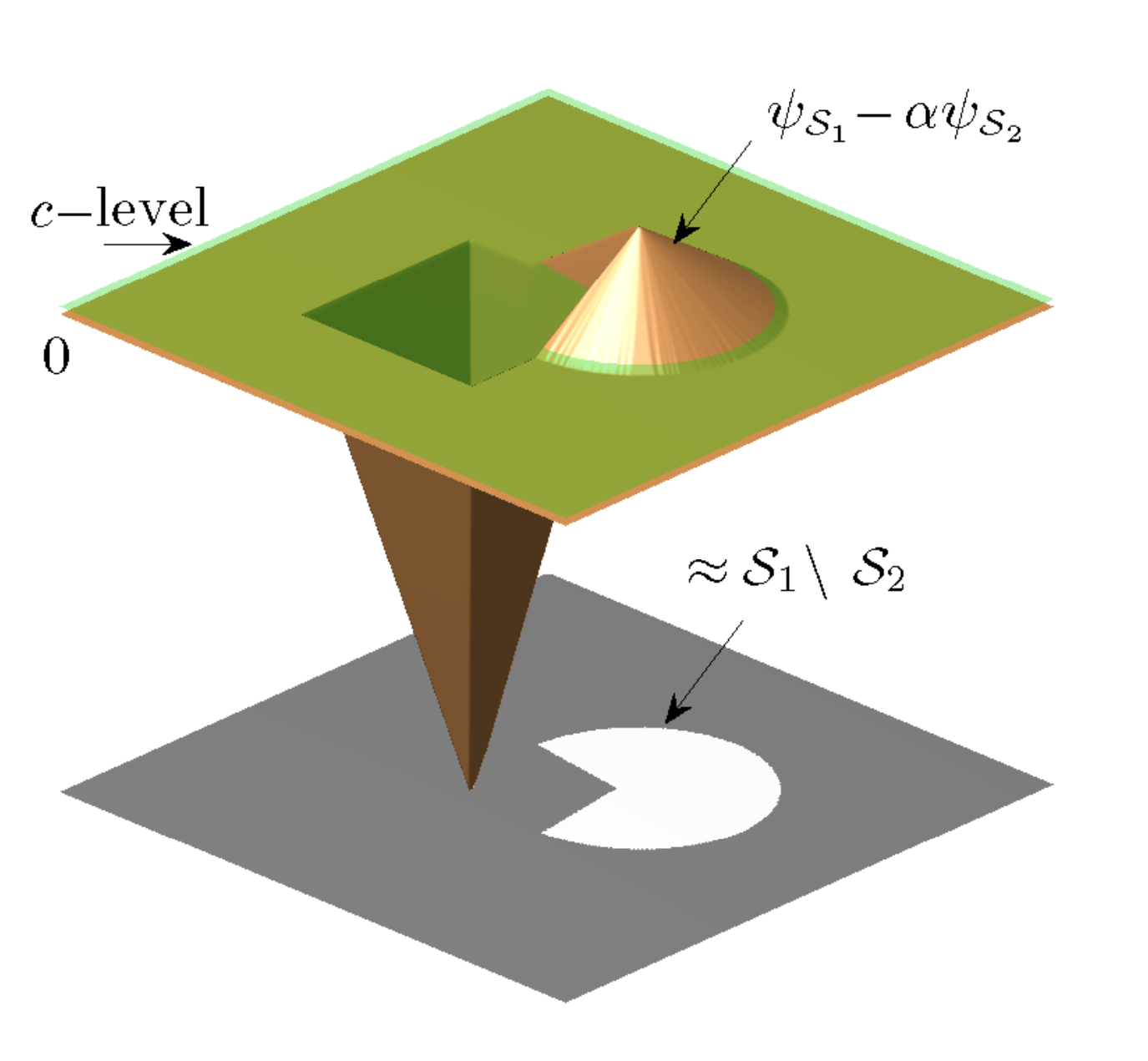}}
\caption{The pseudo-logical behavior of two knolls $\psi_{\mathcal{S}_1}$ and $\psi_{\mathcal{S}_2}$ by considering a close to zero level set $c$ (a) The $c$-level set of $\psi_{\mathcal{S}_1}+\psi_{\mathcal{S}_2}$ approximately represents ${\mathcal{S}_1}\cup{\mathcal{S}_2}$ (b) When $\alpha\gg 1$ the $c$-level set of $\psi_1-\alpha\psi_2$ approximately represents ${\mathcal{S}_1}\setminus{\mathcal{S}_2}$} %
\label{fig2}%
\end{figure}

The idea may be easily incorporated with the notation of level sets by momentarily employing a $c>0$ level set instead of the zero level set. This unusual lifting is due to the compact support of the knolls, as using a zero level set causes ambiguity in identifying the underlying shape boundaries. In this context, for small values $c> 0$, the $c$-level set of $\psi_{\mathcal{S}_1}+\psi_{\mathcal{S}_2}$ approximately represents $\mathcal{S}_1\cup \mathcal{S}_2$ and for large weights $\alpha\gg 1$, the $c$-level set of $\psi_{\mathcal{S}_1}-\alpha\psi_{\mathcal{S}_2}$ approximately represents $\mathcal{S}_1\setminus \mathcal{S}_2$ (see Fig \ref{fig2}). Reverting to (\ref{eq9}), a more general corollary states that for a parametric shape defined as
\begin{equation}\label{eq14}
 \Omega^{\boldsymbol{\alpha}}_{\mathcal{I}_\oplus,\mathcal{I}_\ominus}\triangleq \mbox{supp}^+ \Big( \sum_{i\in \mathcal{I}_\oplus} \alpha_i \psi_{\mathcal{S}_i}(x)-\sum_{j\in \mathcal{I}_\ominus} \alpha_j \psi_{\mathcal{S}_j}(x)-c\Big),
\end{equation}
making
\begin{equation}\label{eq14a}
{c}/{\alpha_i} \to 0^+ \; \mbox{and} \; {\alpha_j}/{\alpha_i} \to +\infty, \qquad (\mbox{for all}\; i\in \mathcal{I}_\oplus, j\in \mathcal{I}_\ominus),
\end{equation}
results
\begin{equation}\label{eq16}
 \mathscr{D}( \Omega^{\boldsymbol{\alpha}}_{\mathcal{I}_\oplus,\mathcal{I}_\ominus},  \Omega_{\mathcal{I}_\oplus,\mathcal{I}_\ominus}) \to 0.
\end{equation}
In other words, by considering the $c$-level set, the parametric function $\sum_{i=1}^{n_d}\alpha_i \psi_{\mathcal{S}_i}(x)$ is capable of producing shapes arbitrarily close to $\Omega_{\mathcal{I}_\oplus,\mathcal{I}_\ominus}$. To conform with prevalent level set methods we henceforth consider the zero level set by transferring the lifting parameter $c$ into the level set formulation as
\begin{equation}\label{eq17}
\phi(x,\boldsymbol{\alpha})=-c+\sum_{i=1}^{n_d}\alpha_i \psi_{\mathcal{S}_i}(x).
\end{equation}
Of course combining the shapes by varying the coefficients $\boldsymbol{\alpha}$ in (\ref{eq17}) brings continuum into the problem. In other words, instead of solving a combinatorial problem and searching among a certain number of possibilities, our search is performed in a subspace spanned by $\psi_{\mathcal{S}_i}(x)$, which contains elements that can represent shapes arbitrarily close to $\Omega_{\mathcal{I}_\oplus,\mathcal{I}_\ominus}$. Searching in this subspace not only relaxes the combinatorial problem and makes it capable of applying descent search methods, also provides the option of exploring shapes that are not among $\Omega_{\mathcal{I}_\oplus,\mathcal{I}_\ominus}$ possibilities (Fig \ref{fig4}).

\subsection{Generating Knolls from Given Shapes} 
\begin{figure}[!htbp]
\begin{minipage}[b]{0.47\linewidth}\hspace{-.4cm}
\includegraphics[width=64mm]{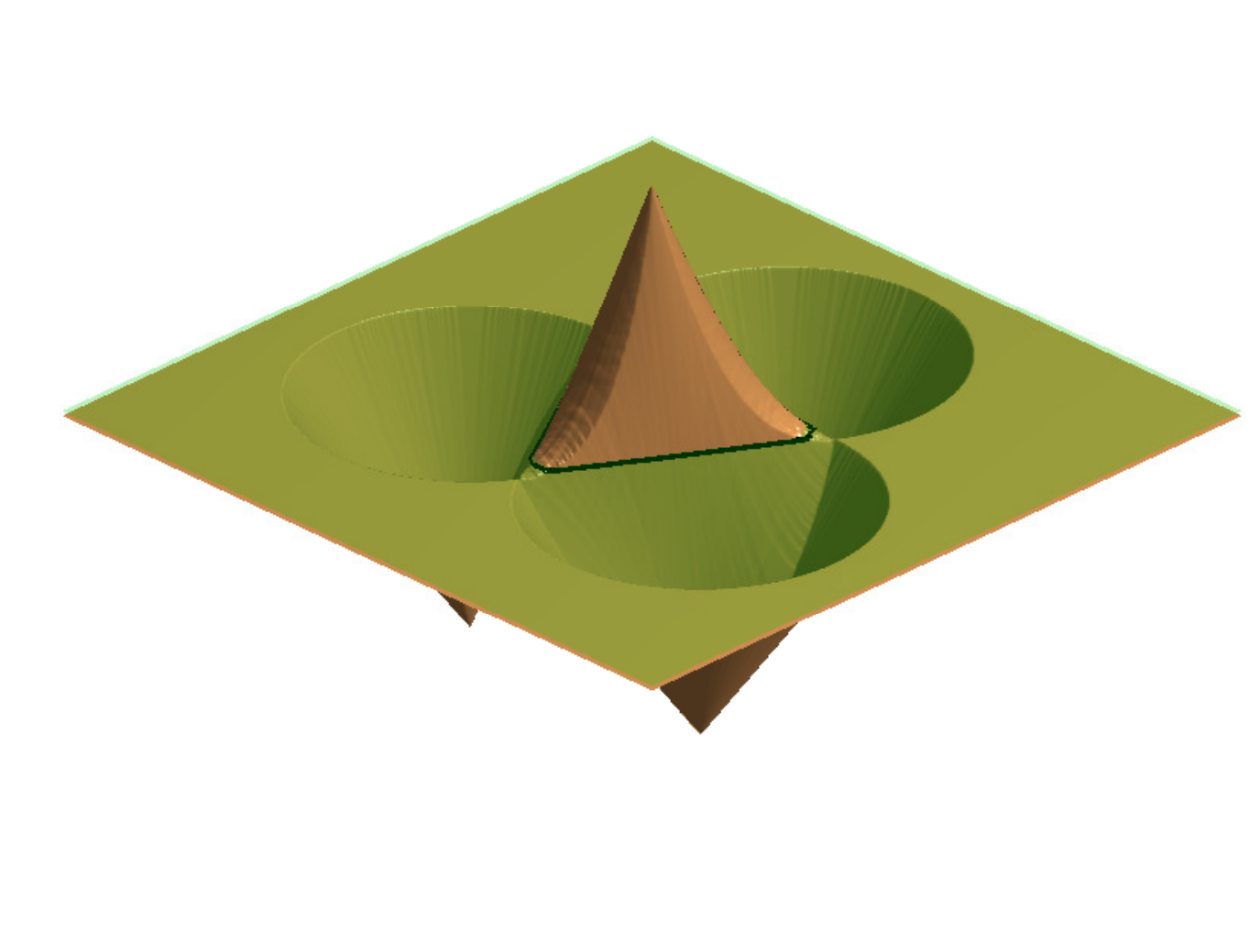}
\caption{Interaction of four identical circular knolls with different signs which gives rise to a triangular region with almost straight sides. The resulting shape is certainly not in $\Omega_{\mathcal{I}_\oplus,\mathcal{I}_\ominus}$ possibilities}
\label{fig4}
\end{minipage}
\hspace{4mm}
\begin{minipage}[b]{0.47\linewidth}
\includegraphics[width=64mm]{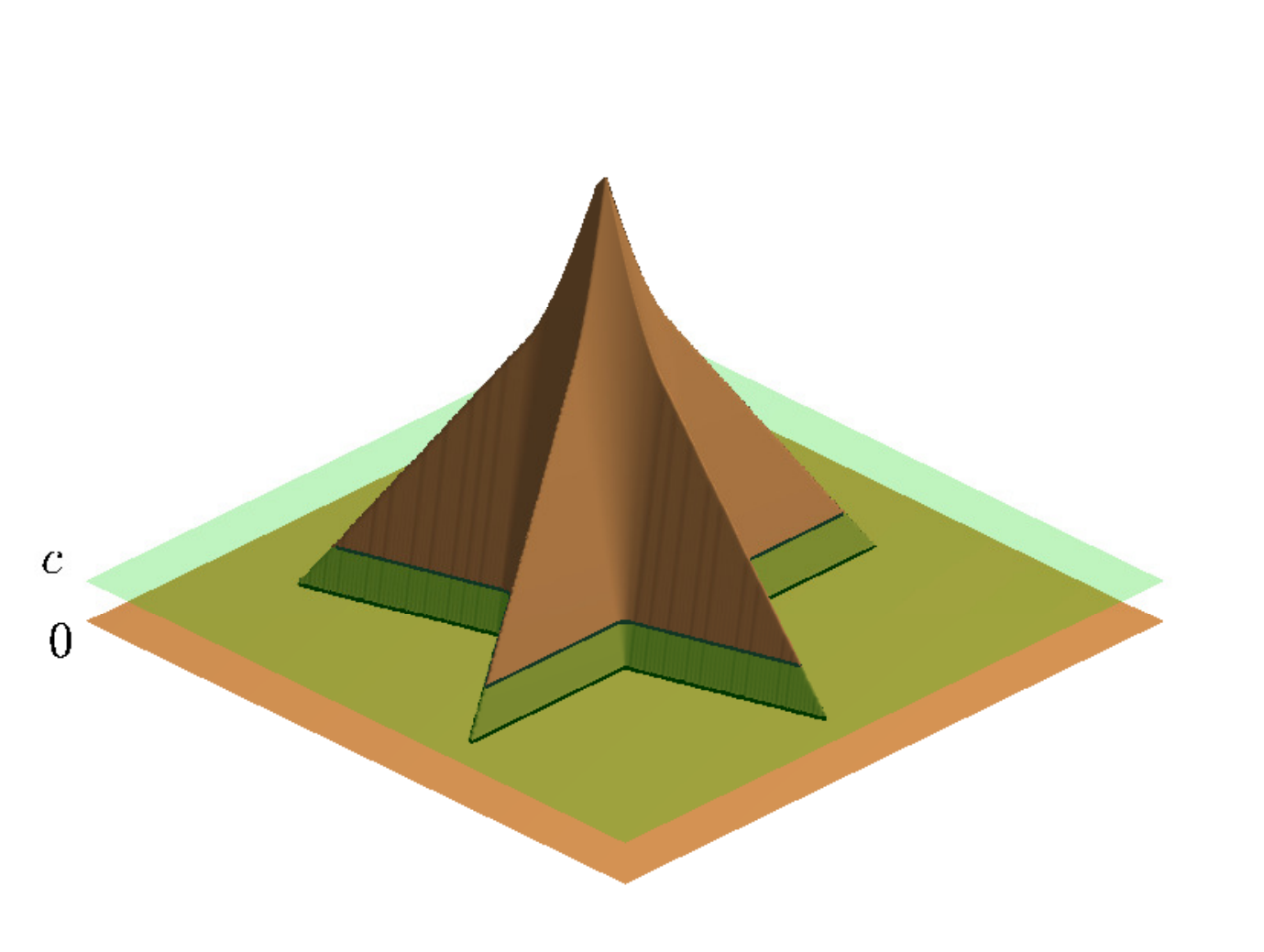}
\caption{A knoll corresponding to a star-shaped region. The shape resulted at some $c$-level set near zero inherits many geometrical features of the knoll's support such a the general structure and corners}
\label{fig3}
\end{minipage}
\end{figure}

Earlier in Section 1.1, we pointed out the use of signed distance functions (SDF) in re-initialization of traditional level sets. In this context at every time step of the Hamilton-Jacobi equation the level set function is re-initialized as an SDF of the form
\begin{equation}\label{eq18}
\rho_{\mathcal{S}}(x)= \left\{
     \begin{array}{cl}
        \mathrm{d}(x,\mathcal{C}) &  x \in \mathcal{S} \\
        -\;\mathrm{d}(x,\mathcal{C}) &  x \in D\setminus \mathcal{S}
     \end{array},
  \right.
\end{equation}
where $\mathcal{S}$ is the underlying shape and $\mathrm{d}(x,\mathcal{C})$ is the distance between the point and the shape boundaries $\mathcal{C}$. A similar idea may be used to uniquely generate knolls from a given shape. More specifically a knoll may be formed as
\begin{equation}\label{eq19}
\psi_{\mathcal{S}}(x)= \left\{
     \begin{array}{cl}
        \mathrm{d}(x,\mathcal{C}) &  x \in \mathcal{S} \\
        0 &  x \in D\setminus \mathcal{S}
     \end{array},
  \right.
\end{equation}
which basically implies $\psi_{\mathcal{S}}(x)=\rho_{\mathcal{S}}(x)^+$. Besides their widespread use in the level set community, our intention of using distance functions is to provide a rather basic expression for the knolls and benefit the fact that shapes at the $c$-level sets of such knolls uniformly inherit many geometrical features from the original shape (see Fig \ref{fig3}). Moreover, since the the $c$-level set of a knoll $\psi_{\mathcal{S}}$ (or equivalently the zero level set of $\psi_{\mathcal{S}}-c$) is slightly smaller than $\mathcal{S}$, in a level set representation as (\ref{eq17}), knolls may be generated with slightly larger supports, e.g., $\psi_{\mathcal{S}}(x)=\big(\rho_{\mathcal{S}}(x)+c\big)^+$. By this choice, reconstructing exact elements of the dictionary does not require having very large weights associated with their knolls.

\section{Sparse Reconstruction}
We so far discussed how considering a parametric level set as (\ref{eq17}) relaxes (\ref{eq10}) into a minimization of the form
\begin{equation}\label{eq20}
\tilde{ \boldsymbol{\alpha}}=\operatorname*{arg\,min}_{\|\balpha\|_0\leq s} \mathcal{E}(\boldsymbol{ \alpha}),
\end{equation}
where $\mathcal{E}(\boldsymbol{ \alpha})=\mathpzc{E}(\Omega_{ \phi(x,\boldsymbol{ \alpha})})$. We are certainly interested in sparse solutions of $\balpha$ to avoid shape redundancy and ill-conditioning associated with large dictionaries. In general $\mathcal{E}$ is a nonlinear and non-convex function of $\boldsymbol{\alpha}$, we however assume it to be a sufficiently smooth function of $\boldsymbol{\alpha}$.

In the sequel we provide a general overview of sparse recovery techniques in finding solutions of underdetermined linear systems. We then use this notion to develop our sparsity promoting algorithm applicable to a large class of functionals appeared in imaging applications.

\subsection{Background on Sparse Recovery Techniques}
Finding sparse solutions of linear systems is a broad area of research in imaging science \cite{romberg2008imaging, goldstein2009split, lustig2007sparse}. For an underdetermined linear system, $\boldsymbol{A}\balpha=\boldsymbol{b}$, the main underlying problem is
\begin{equation}\label{eq21}
\mbox{minimize} \quad \|\balpha\|_0 \qquad \mbox{s.t.} \qquad \boldsymbol{A}\balpha=\boldsymbol{b}.
\end{equation}
The matrix $\boldsymbol{A}$ is $m$-by-$n$ where $m\ll n$ and $\boldsymbol{b}$ is a vector of length $m$. Problem (\ref{eq21}) is in general a hard combinatorial problem. However, it was brought into attention that relaxing (\ref{eq21}) by replacing $\|\balpha\|_0$ with $\|\balpha\|_1$, known as the \emph{basis pursuit} (BP) problem \cite{chen2001atomic}, can still result in sparse solutions. Moreover, under certain conditions on $\boldsymbol{A}$, a BP solution perfectly coincides with the solution of (\ref{eq21}) \cite{candes2006robust, donoho2006compressed, candes2006stable}. In case of noisy observations, the linear equality in (\ref{eq21}) is replaced with a least squares inequality as
\begin{equation}\label{eq22}
\mbox{minimize} \quad \|\balpha\|_1 \qquad \mbox{s.t.} \qquad \|\boldsymbol{A}\balpha-\boldsymbol{b}\|_2\leq \sigma,
\end{equation}
where $\sigma$ is the noise level. This convex minimization is known as the basis pursuit denoising (BPDN) problem. We consider the BP problem as a specific case of BPDN when $\sigma=0$.

Another convex problem earlier used to induce sparsity on the solutions of a linear system is the \emph{least absolute shrinkage and selection operator} (Lasso) \cite{tibshirani1996regression}, phrased as
\begin{equation}\label{eq23}
\mbox{minimize} \quad \|\boldsymbol{A}\balpha-\boldsymbol{b}\|_2 \qquad \mbox{s.t.} \qquad \|\balpha\|_1 \leq \tau.
\end{equation}
For a certain value $\tau=\tau_\sigma$, problem (\ref{eq23}) becomes equivalent to (\ref{eq22}), however, in general $\tau_\sigma$ cannot be determined a priori. To promote sparsity, in most applications solving (\ref{eq22}) is more desirable than solving (\ref{eq23}). This is mainly due to the fact that a good estimate of the noise in the measurements is more likely to be available than prior knowledge of $\tau$, the $\ell_1$-magnitude of a sparse solution.

A variety of solution strategies may be taken for both the BPDN and Lasso problems (e.g., see \cite{osborne2000new, candes2005l1, figueiredo2007gradient, tibshirani1996regression,osborne2000lasso}). Relatively speaking, the smooth cost associated with the Lasso makes it an easier problem to approach compared to the BPDN which involves a non-smooth cost minimization. Minimization of smooth costs over convex sets may be efficiently handled using methods such as spectral projected gradient (SPG) \cite{birgin2000nonmonotone}.

We briefly overview the SPG-based technique proposed by van den Berg and Friedlander in \cite{van2008probing}, who show that a numerically efficient way of solving the BPDN problem is by solving a series of Lasso problems. This technique plays a major role in developing the algorithm that we put forth later in this paper. Their method also allows use of $\boldsymbol{A}$ as a linear operator with no explicit matrix form which suits our general representation.

\begin{figure}[!htbp]
\centering \includegraphics[width=85mm]{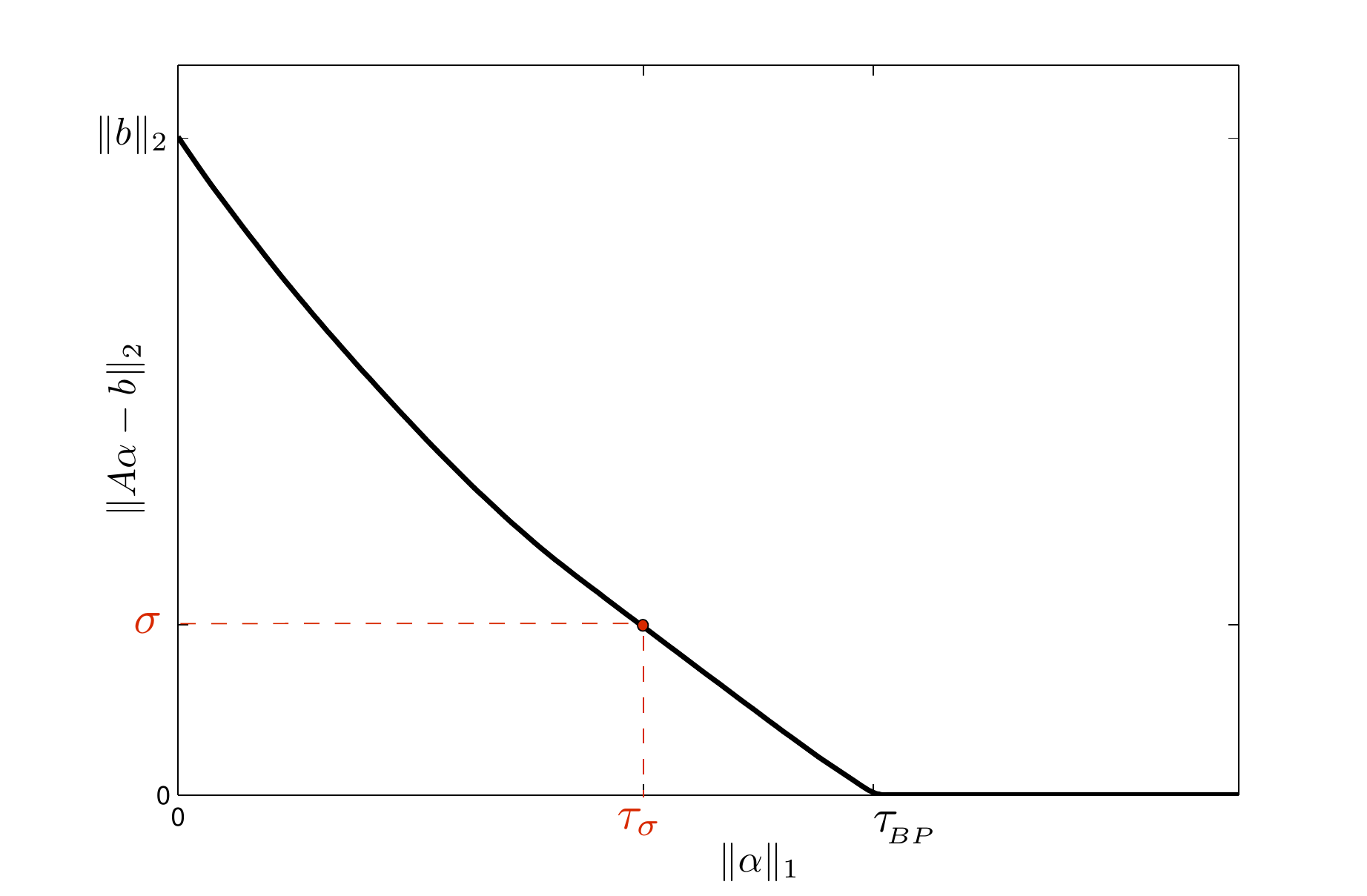}
\caption{A typical Pareto curve corresponding to the Lasso and BPDN problems. For a given $\sigma$, the corresponding value $\tau_\sigma$ makes Lasso and BPDN problems share solution. For values of $\tau$ larger than $\tau_{BP}$ (the $\ell_1$-norm of the BP solution), the residual $\|A\alpha-b\|_2$ is zero.}%
\label{fig5}%
\end{figure}

Following \cite{van2008probing}, for a given $\tau\geq 0$, a single parameter function $\varphi(\tau)$ is defined as
\begin{equation}\label{eq24}
\varphi(\tau)\triangleq\|\boldsymbol{b}-\boldsymbol{A}\balpha_\tau\|_2,
\end{equation}
where $\balpha_\tau$ is a solution of (\ref{eq23}) for the given $\tau$. It should be noted that multiple solutions to the Lasso may exist, however, problem's convexity requires all of them to generate an equal error term $\varphi(\tau)$. It is shown in \cite{van2008probing} that $\varphi(.)$ is a nonincreasing convex function which is continuously differentiable and
\begin{equation}\label{eq25}
\varphi'(\tau)=\frac{\mbox{d}\varphi}{\mbox{d}\tau}=-\frac{\|\boldsymbol{A}^T(\boldsymbol{b}-\boldsymbol{A}\balpha_\tau) \|_\infty}{\|\boldsymbol{b}-\boldsymbol{A}\balpha_\tau \|_2}.
\end{equation}
The graph of $\varphi(\tau)$ in terms of $\tau$ provides an optimal trade-off between the residual error $\|\boldsymbol{b}-\boldsymbol{A}\balpha_\tau \|_2$ and the $\ell_1$-norm of the solution, forming the Pareto curve (Fig \ref{fig5}). The Lasso and BPDN problems are basically two different characterizations of this curve. In other words, for a given $\sigma$ a solution to
\begin{equation}\label{eq26}
\varphi(\tau)=\sigma,
\end{equation}
returns the value of $\tau_\sigma$ that makes solutions of BPDN and Lasso problems coincide. As a matter of fact solving (\ref{eq26}) does not require full access to the Pareto curve. The convex nature of the problem and a closed form expression for $\varphi'(\tau)$ are sufficient to employ root finding techniques and acquire $\tau_\sigma$. In this context, a Newton iterative process $\tau_{\ell+1}=\tau_\ell+\Delta\tau_\ell$, where
\begin{equation}\label{eq27}
\Delta\tau_\ell=\frac{\sigma- \varphi(\tau_\ell)}{\varphi'(\tau_\ell)}, \qquad \ell=0,1,\cdots
\end{equation}
is capable of generating a sequence of Lasso parameters $\tau_\ell$ that superlinearly converge to $\tau_\sigma$. This result makes Lasso the central tool in solving the BPDN problem. More specifically, for a given $\sigma\geq 0$, a solution to the BPDN problem is acquired by solving a series of Lasso problems parameterized by $\tau_\ell$ and ultimately arriving at a Lasso that is parameterized by $\tau_\sigma$, which basically solves the BPDN problem. To solve the underlying Lasso problems efficiently authors in \cite{van2008probing} employ the SPG technique detailed in \cite{birgin2000nonmonotone}.

To address a broader class of problems, specifically in a variational framework, in our future formulations we look beyond matrix representatives of $\boldsymbol{A}$. We consider a more general case that $A: \mathbb{R}^n \to \mathbb{S}$ is a bounded linear operator and $\mathbb{S}$ represents a real Hilbert space\footnote{We therefore keep the notation general and do not use a bold syntax (representative of a matrix form) for the linear operator and elements of $\mathbb{S}$.}. An example of $\mathbb{S}$, other than $\mathbb{R}^m$, is the space of bounded continuous functions $C(\mathbb{R})$. This generalization only requires slight modifications to the above overview of sparse recovery techniques, basically by replacing $\|.\|_2$ with $\|.\|_\mathbb{S}$, the norm induced by the Hilbert space $\mathbb{S}$, and using the adjoint operator $A^*$ in place of $\boldsymbol{A}^T$. The sparsity constraints applied to $\balpha$ remain intact as the domain of $A$ is still $\mathbb{R}^n$. For linear systems, the technique proposed in \cite{van2008probing} makes such generalization possible and conveniently handles cases that $A$ is not explicitly available as a matrix.

\subsection{A More General Sparsity Constrained Optimization}\label{sec-Im}
 The majority of sparsity-constrained optimization techniques are merely applicable to linear systems and corresponding quadratic costs. We are however interested in sparse stationary points of $\mathcal{E}({\balpha})$ in (\ref{eq20}), which is not generally a quadratic function of $\balpha$. This is a more complex problem and a new area of research which still requires further development. The current available techniques are mainly in the form of (\ref{eq20}) for which we are required to know the degree of sparsity, $s$, prior to the minimization \cite{beck2012sparsity, blumensath2012compressed, bahmani2011greedy}.

 For the purpose of this paper, we consider a particular form of $\mathcal{E}({\balpha})$ applicable to a large class of imaging problems. We then provide a sparsity promoting minimization scheme that takes a rather short iterative process to converge and does not require prior assignment of $s$. Accordingly, for a given real Hilbert space $\mathbb{S}$, we consider $\mathcal{E}(\balpha)$ to take a form (or a finite sum of forms) as
\begin{equation}\label{eq28}
\mathcal{E}(\balpha)=\|\mathcal{G}(\balpha)\|_\mathbb{S}^2= \big \langle \mathcal{G}(\balpha), \mathcal{G}(\balpha)\big \rangle_\mathbb{S},
\end{equation}
where $\mathcal{G}: \mathbb{R}^n\to \mathbb{S}$ is a Fr\'{e}chet differentiable map applied to $\balpha$. The classic compressed sensing problem can be interpreted as a special case of (\ref{eq28}) when $\mathcal{G}(.)$ is an affine function of $\balpha$, i.e., $\mathcal{G}(\balpha)=A\balpha-b$. In the following we present two classes of imaging applications that follow the form (\ref{eq28}) and will be later considered in Section \ref{exsec}.

\subsubsection{Shape-Based Inverse Problems} For this class of problems, in the basic case, the spatial parameter of interest $u$ is modelled as
\begin{equation}\label{eq29}
u(x)=\left\{
     \begin{array}{cl}
        u_{in} &  x \in \Omega \\
        u_{ex} &  x \in D\setminus \Omega
     \end{array},
  \right.
\end{equation}
where $u_{in}$ corresponds to the texture value of an inclusion (the shape) and $u_{ex}$ represents the background texture. The texture values may be represented by scalars or low order parametric models which may/may not be known a priori \cite{aghasi2011parametric}. Using a parametric function $\phi(x,\balpha)$ to characterize the shape by its zero level set we may rewrite the parameter of interest as
\begin{equation}\label{eq30}
u(x,\balpha)=u_{in} H\big(\phi(x,\balpha)\big)+u_{ex} \Big(1-H\big(\phi(x,\balpha)\big)\Big),
\end{equation}
where $H(.)$ is the heaviside function, usually replaced with a smooth approximation $H_{rg}(.)$ to maintain differentiability \cite{aghasi2011parametric}. The inverse problem may now be formulated as minimizing the model-data mismatch to acquire the shape parameters $\balpha$, i.e.,
\begin{equation}\label{eq31}
\operatorname*{min}_{\balpha} \big \|v-\mathcal{M}\big(u(x,\balpha)\big)\big \|_{\mathbb{S}_v}^2.
\end{equation}
For a smooth physical model $\mathcal{M}$ that relates the parameter of interest $u$ to the observations, the underlying inversion cost may be interpreted as a functional of the form (\ref{eq28}) when $\mathcal{G}(\balpha)=v-\mathcal{M}\big(u(x,\balpha)\big)$.

\subsubsection{Image Segmentation} \label{sec-Im2}  A well known variation form to segment an image $u(x)$ into two disjoint regions $\Omega$ and $D\setminus \Omega$ corresponds to a functional
\begin{equation}\label{eq32}
\mathpzc{E}(\Omega)=\mathpzc{E}_r(\Omega)+\int_{\Omega}r_{in}\big(u(x)\big)\;\mbox{d}x + \int_{D\setminus \Omega}r_{ex}\big(u(x)\big)\;\mbox{d}x,
\end{equation}
where $\mathpzc{E}_r(\Omega)$ represents some regularity constraints on the segmented regions (e.g., smoothness and compactness \cite{caselles1997geodesic, chan2001active}) and $r_{in}(.)\geq 0$ and $r_{ex}(.) \geq 0$ are some inhomogeneity measures of $u$ in each region. A well known instance of such measures is observed in the work by Chan and Vese \cite{chan2001active}, which suggests $r_{in}(u)=(u(x)-\tilde u_{in})^2$ and $r_{ex}(u)=(u(x)-\tilde u_{ex})^2$. Here $\tilde u_{in}$ and $\tilde u_{ex}$ are scalars representing the average pixel value within each region. Another example is using a maximum likelihood approach for the pixel intensities to arrive at
\begin{equation}\label{eq33}
r_{\chi}(u)=-\log \Big( p_{\chi}\big(u(x)\big) \Big), \qquad \chi:in,ex,
\end{equation}
where $p_{in}(u)$ and $p_{ex}(u)$ are pixel intensity distributions \cite{paragios2002geodesic, cremers2007review}. The scalar parameters $\tilde u_{in}$, $\tilde u_{ex}$ and the distributions $p_{in}(.)$, $p_{ex}(.)$ may be known a priori or may be estimated in the course of segmentation. Considering the image dependent terms in (\ref{eq32}), a parametric level set function may be employed to form the segmentation functional
\begin{align}\nonumber
\mathcal{E}(\balpha)=&\int_{D}r_{in}\big(u(x)\big)H\big(\phi(x,\balpha)\big)\;\mbox{d}x \\& + \int_{D}r_{ex}\big(u(x)\big)\Big(1- H\big(\phi(x,\balpha)\big) \Big)\;\mbox{d}x.\label{eq34}
\end{align}

We note that in using a parametric level set technique the regularity constraints may either be neglected thanks to the intrinsic smoothness of the parametric function \cite{gelas2007compactly}, or in our case taken into consideration via the knolls sparsity. The energy functional in (\ref{eq34}) may be written as sum of two functionals $\mathcal{E}_{in}(\balpha)$ and $\mathcal{E}_{ex}(\balpha)$ in form of (\ref{eq28}) using $\mathcal{G}_{in}(\balpha)= \sqrt{r_{in}(u)H(\phi(x,\balpha))}$ and $\mathcal{G}_{ex}(\balpha)=\sqrt{r_{ex}(u)(1- H(\phi(x,\balpha)) )}$. Employing a smooth function $H_{rg}(.):\mathbb{R}\to [0,1]$ to approximate the Heaviside function can guarantee smoothness of $\mathcal{G}_{in}(.)$ and $\mathcal{G}_{ex}(.)$ in $\balpha$.

\subsection{A Sparse Nonlinear Minimization Technique}
To minimize a sufficiently smooth cost function $\mathcal{E}(\balpha):\mathbb{R}^n\to \mathbb{R}$, a well known iterative scheme is Newton's method. Starting from an initial vector $\balpha_0$, Newton's method proceeds by generating $\balpha_k$ vectors ($k=1,2,\cdots$) that progressively reduce the cost to reach a minima. Having $\balpha_k$ available, the successive vector is written as $\balpha_{k+1}=\balpha_k+\bdelta_k$ (or a multiple of $\bdelta_k$), where the step is determined by minimizing the second order Taylor approximation of the cost around $\balpha_k$:
\begin{equation}\label{eq35}
\bdelta_k= \operatorname*{arg\,min}_{\bdelta}\; \mathcal{E}(\balpha_k)+ \je(\balpha_k)\bdelta +\frac{1}{2}\bdelta^T \he(\balpha_k)\bdelta.
\end{equation}
Here $\je(\balpha_k)$ denotes the Jacobian vector at $\balpha_k$ and $\he(\balpha_k)$ is the corresponding Hessian matrix.

A positive definite Hessian guarantees existence of a minima for (\ref{eq35}) which leads to the closed form expression $\bdelta_k=( \he(\balpha_k))^{-1}\je(\balpha_k)$. Most Newton type techniques either impose such positivity constraint on the Hessian or approximate it with an at least positive semi-definite matrix to guarantee convexity of the underlying cost \cite{dennis1996numerical}.

To consider applying a Newton type technique to the energy functional in (\ref{eq28}), using calculus of operators \cite{griffel2002applied,aghasi2011parametric}, we have
\begin{equation}\label{eq36}
 \je(\balpha)\bdelta = 2\;\big \langle \mathcal{G}'(\balpha)\bdelta, \mathcal{G}(\balpha)\big \rangle_\mathbb{S},
\end{equation}
and
\begin{equation}\label{eq37}
\frac{1}{2}\bdelta^T\he(\balpha)\bdelta = \big \langle \mathcal{G}'(\balpha)\bdelta, \mathcal{G}'(\balpha)\bdelta\big \rangle_\mathbb{S}+\big \langle \mathcal{G}''(\balpha)[\bdelta,\bdelta], \mathcal{G}(\balpha)\big \rangle_\mathbb{S}.
\end{equation}
Here $\mathcal{G}'(\balpha)[.]:\mathbb{R}^n\to \mathbb{S}$ is a linear operator representing the first order Fr\'{e}chet derivative of $\mathcal{G}$ at $\balpha$ and $\mathcal{G}''(\balpha)[.,.]:\mathbb{R}^n\times \mathbb{R}^n\to \mathbb{S}$ is a bilinear operator that represents the second order Fr\'{e}chet derivative of $\mathcal{G}$. When $\mathbb{S}$ is taken to be $\mathbb{R}^m$, $\mathcal{G}'$ and $\mathcal{G}''$ are simply the Jacobian matrix and the Hessian tensor of $\mathcal{G}$. For a positive semi-definite approximation to $\he$, in Gauss-Newton methods \cite{dennis1996numerical} the inner product term containing $\mathcal{G}''$ is neglected in (\ref{eq37}). Therefore the underlying quadratic cost corresponding to (\ref{eq35}) becomes
\begin{align}\label{eq38}\nonumber
\big \| \mathcal{G}(\balpha_k)\big\|_\mathbb{S}^2 &+ 2\;\big \langle \mathcal{G}'(\balpha_k)\bdelta, \mathcal{G}(\balpha_k)\big \rangle_\mathbb{S} + \big \| \mathcal{G}'(\balpha_k)\bdelta\big \|_\mathbb{S}^2\\&=\| \mathcal{G}'(\balpha_k)\bdelta+\mathcal{G}(\balpha_k)\|_\mathbb{S}^2,\nonumber
\end{align}
forming a convex function of $\bdelta$. Although the underlying quadratic cost is convex, the linear operator $\mathcal{G}'$ may be ill conditioned and some type of regularization may be required to obtain $\bdelta_k$. For instance Levenberg-Marquardt algorithm is a variant of Gauss-Newton technique that uses a Tikhonov type regularizer to minimize the underlying quadratic cost \cite{dennis1996numerical}.

We however consider a different regularizer which is capable of promoting sparsity on consecutive estimates of $\balpha_k$. More specifically having $\balpha_k$ available we regularize the subproblem by promoting sparsity on the next potential estimate $\balpha_{k+1}$ through a BPDN problem, i.e.,
\begin{equation}\label{eq39}
\bdelta_k =\operatorname*{arg\,min}_{\bdelta} \quad \|\balpha_k+\bdelta\|_1 \qquad \mbox{s.t.} \qquad \| \mathcal{G}'(\balpha_k)\bdelta+\mathcal{G}(\balpha_k)\|_\mathbb{S}\leq \sigma.
\end{equation}
This is certainly a translated version of the BPDN problem and a simple change of variable as $\boldsymbol{\theta} =\balpha_k+\bdelta$ converts it to the standard form (\ref{eq22}). As $\bdelta$ is expected to be very small about an accumulation point, a reasonable choice of $\sigma$ is an estimate of $\|\mathcal{G}(\balpha)\|_\mathbb{S}^2$ near a minima. In Appendix A we show that when $\sigma=\min_{\balpha} \|\mathcal{G}(\balpha)\|_\mathbb{S}^2$, or even more locally when $\sigma\leq \|\mathcal{G}(\balpha_k)\|_\mathbb{S}^2$, the direction $\bdelta_k$ acquired in (\ref{eq39}) is a descent direction.

In the case of inverse problems, $\sigma$ may be an estimate of noise power in the measurements. Nevertheless, if such estimate is not available, a substitute for (\ref{eq39}) would be a parameter free BP problem
\begin{eqnarray}\nonumber
&\bdelta_k =\operatorname*{arg\,min}_{\bdelta} \quad \|\balpha_k+\bdelta\|_1\\ &\mbox{s.t.} \qquad  \mathcal{G}'^*(\balpha_k)\mathcal{G}'(\balpha_k)\bdelta+\mathcal{G}'^*(\balpha_k)\mathcal{G}(\balpha_k)=0,\label{eq40}
\end{eqnarray}
corresponding to the classic Gauss-Newton equation with an additional sparsity constraint. If $\mathcal{G}'(\balpha_k)\bdelta+\mathcal{G}(\balpha_k)=0$ has a class of solutions, which is very likely for an ill-conditioned linear system, this equality constraint may be replaced with the one in (\ref{eq40}) and make it equivalent to the case $\sigma=0$ in (\ref{eq39}).

As stated in Section 3.1, at every iteration a minimizer to (\ref{eq39}) may be determined via a sequence of Lasso solves
\begin{eqnarray}\label{eq41}\nonumber
 &\bdelta_{k,\ell} =\operatorname*{arg\,min}_{\bdelta} \quad \| \mathcal{G}'(\balpha_k)\bdelta+\mathcal{G}(\balpha_k)\|_\mathbb{S}\\ &\mbox{s.t.} \qquad \|\balpha_k+\bdelta\|_1\leq \tau_{k,\ell},
\end{eqnarray}
where for $\ell=0,1,2,\cdots$ the $\ell_1$ radius $\tau_{k,\ell}$ is updated according to (\ref{eq27}) and sequentially $\lim_{\ell\to \infty}\bdelta_{k,\ell}=\bdelta_{k}$. Once the $k$-th iteration is complete, a similar Lasso search is performed along the Pareto curve corresponding to $\| \mathcal{G}'(\balpha_{k+1})\bdelta+\mathcal{G}(\balpha_{k+1})\|_\mathbb{S}$ to find $\bdelta_{k+1}$.

In determining the consecutive directions $\bdelta_k$, it is computationally desirable to perform incomplete iterations on $\ell$ and make a more accurate search for $\bdelta_k$ when $\balpha_k$ approaches a sparse solution. For this purpose, it is also shown in Appendix A that if for some $\ell$, $\tau_{k,\ell}\geq \|\balpha_k\|_1$, the corresponding step, $\bdelta_{k,\ell}$, is a descent direction and can be used as a step to reduce the nonlinear cost.

Also, for small values of $k$ that $\balpha_k$ may be far from a sparse solution, $\tau_{k,\ell}$ may tend to very large values because of the high cardinality of $\balpha_k$. To prevent this, the values $\tau_{k,\ell}$ may be controlled by a loose threshold value $\tau_{mx}$ and the iterations on $\ell$ may be broken if this threshold is reached. The descent property of the incomplete step for this case is also inferred by a similar argument as the one just mentioned. These two strategies provide us with an indication of when we are eligible to break the $\ell$-iterations and still be certain that the descent property of resulting direction is maintained.

Algorithm 1 provides a detailed picture of the proposed minimization scheme where a line search technique is employed to warrant a proper convergence behavior. We specifically use Armijo rule to determine the step size, which requires the reduction in the cost to be sufficiently large \cite{bertsekas1999nonlinear}. Parameters $\gamma$ and $\beta$ correspond to this reduction and usually take values as $\gamma\in [10^{-5},0.1]$ and $\beta\in [0.1,0.5]$, detailed in \cite{bertsekas1999nonlinear}. The values $\varepsilon_1$ and $\varepsilon_2$ control the size of updates in each block and are very small numbers in the order of machine precision. The value of $\tau_{mx}$ may be taken as a loose overestimate of the $\ell_1$-norm of the solution. Finally $\balpha_0$ is the initialization of unknown parameters and $\tau_0$ is an initial $\ell_1$-ball radius which may arbitrarily taken to be zero or very small.

\begin{algorithm}\label{alg1}
\caption{A sparsity promoting Gauss-Newton algorithm}
\begin{algorithmic}
\STATE {$\mbox{\textbf{input}}\;\;\sigma$};
\STATE \textbf{set} $\gamma$, $\beta$, $\varepsilon_1$, $\varepsilon_2$ and $\tau_{mx}$;
\STATE $\balpha:=\balpha_0$;
\STATE $\tau:=\tau_0$;
\STATE $found:=\mbox{\textbf{false}}$;
\WHILE{$\sim found$}
\STATE $stepfound:=\mbox{\textbf{false}}$;
\WHILE{$\sim stepfound\vspace{.1cm}$}
\STATE{$\bdelta :=\operatorname*{arg\,min}_{\boldsymbol{\eta}} \;\; \| \mathcal{G}'(\balpha)\boldsymbol{\eta}+\mathcal{G}(\balpha)\|_\mathbb{S} \;\; \mbox{s.t.} \;\; \|\balpha+\boldsymbol{\eta}\|_1\leq \tau \vspace{.1cm}$};
\STATE $\phi_\tau'\!\!:=\!\! - \big\| \mathcal{G}'^*(\balpha)[\mathcal{G}'(\balpha)\bdelta+\mathcal{G}(\balpha)]\big\|_\infty  /\| \mathcal{G}'(\balpha)\bdelta+\mathcal{G}(\balpha)\|_\mathbb{S}\vspace{.1cm}$;
\STATE $\Delta\tau:=\big(\sigma- \| \mathcal{G}'(\balpha)\bdelta+\mathcal{G}(\balpha)\|_\mathbb{S}\big)/\phi_\tau'\vspace{.1cm}$;
\STATE $\tau:=\tau+\Delta\tau$;\;\;(alternatively, $\tau:=\min(\tau+\Delta\tau,\tau_{mx})$;)
\STATE $stepfound:=(\tau\geq \|\balpha\|_1) \vee (\Delta\tau\leq \varepsilon_1)$;
\ENDWHILE
\WHILE{$\mathcal{E}(\balpha)-\mathcal{E}(\balpha+\bdelta)< - \gamma \je(\balpha)\bdelta $}
\STATE $\delta:=\beta*\delta$;
\ENDWHILE
\STATE $\balpha:=\balpha+\bdelta$;
\STATE $found:=(\|\bdelta\|_2\leq \varepsilon_2)$;
\ENDWHILE

\end{algorithmic}
\end{algorithm}

\subsection{Employing an Asymmetric Norm}
In the sparsity promoting algorithm proposed we used the $\ell_1$ norm as a convex relaxation to the $\ell_0$ norm. This is a common technique in compressed sensing promoting sparsity on the reconstructions thanks to the sharp vertices of the $\ell_1$ ball. However, in general an $\ell_0$ minimization only targets the cardinality of a vector while an $\ell_1$ minimization also takes into account the vector component values.

\begin{figure*}[t]\hspace{-.5cm}
\subfigure[][]{\includegraphics[width=70mm]{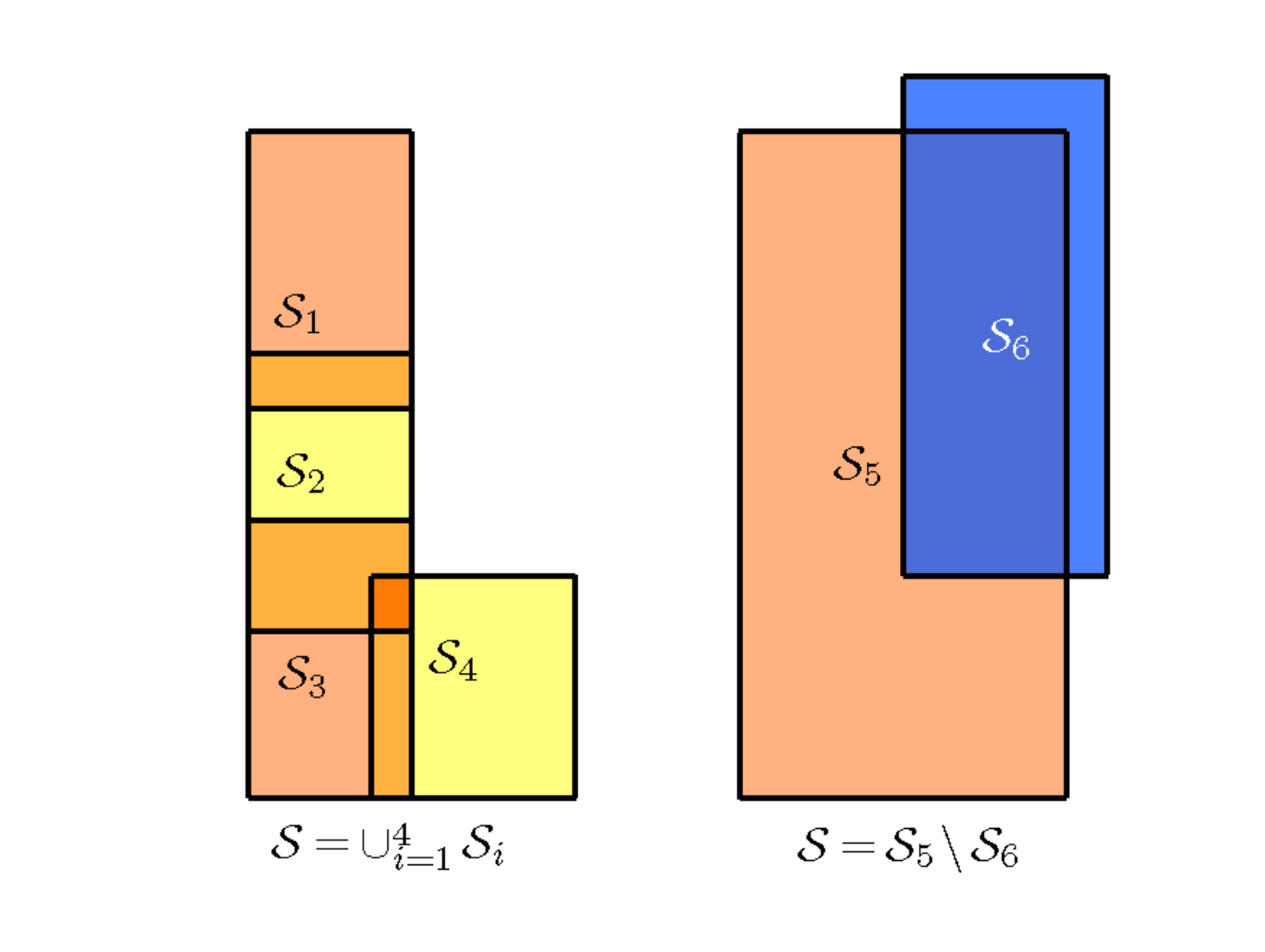}}
\subfigure[][]{\includegraphics[width=65mm]{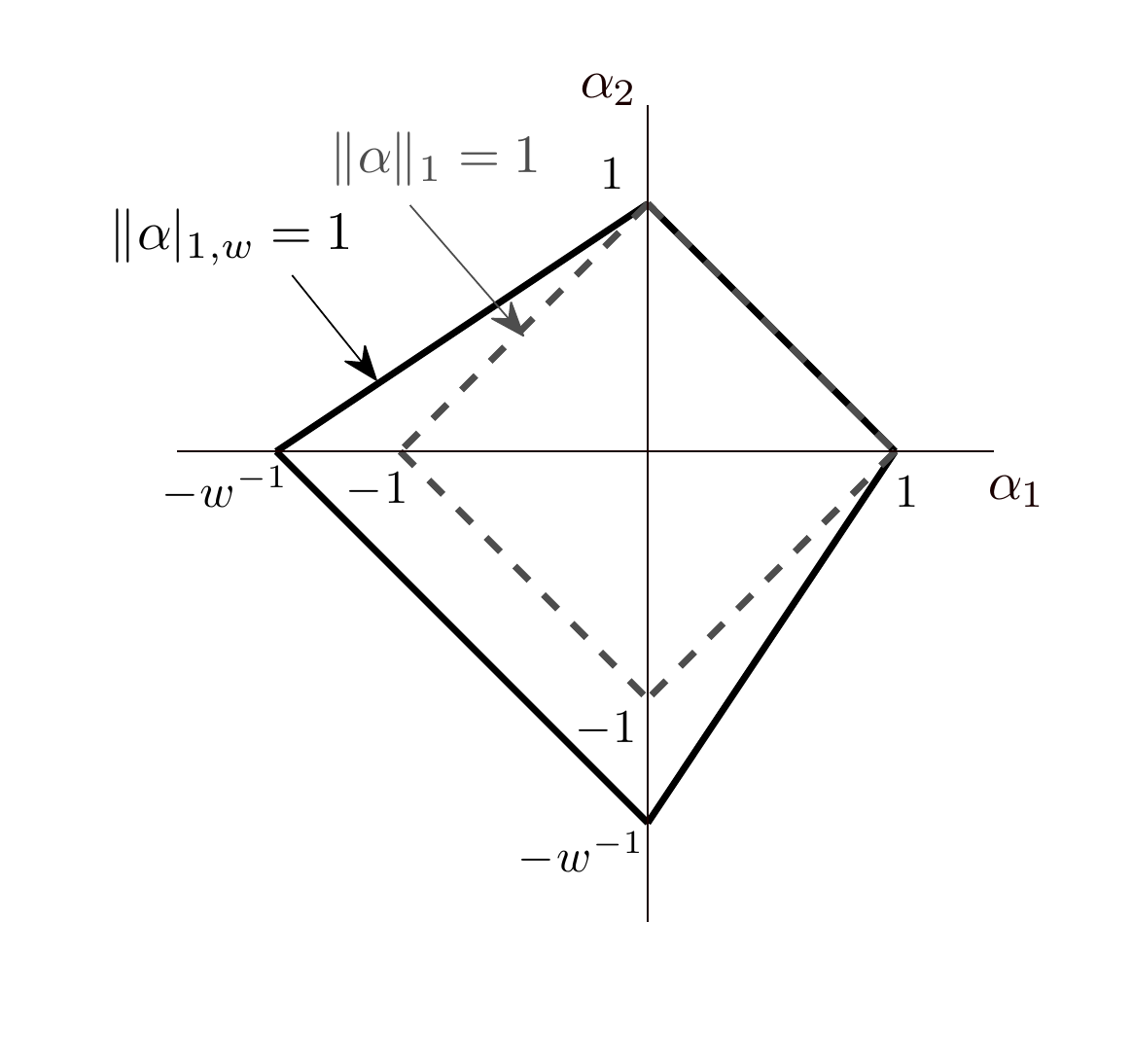}}
\caption{Employing an asymmetric $\ell_1$ norm: (a) Reconstructing an L-shaped region by union of four rectangles (left) and applying the relative complement among two rectangles (right); (b) The $\ell_1$ ball and the asymmetric convex ball corresponding to $\|\balpha|_{1,w}$}\label{fig6}
\end{figure*}

For a parametric representation as (\ref{eq17}), Fig \ref{fig6}(a) shows reconstruction of an L-shaped region where a basic $\ell_1$ minimization does not necessarily provide the best sparse solution. Considering a dictionary that contains all shapes $\mathcal{S}_1,\cdots,\mathcal{S}_6$ shown in the figure, in one case the L-shape is reconstructed via the union of four rectangular knolls while a similar reconstruction is possible by applying the relative complement among only two knolls. When $\|\balpha\|_1$, the $\ell_1$ norm of knoll coefficients, is used to indicate the level of sparsity, the latter case may have a lager $\ell_1$ norm although only two shapes are exploited. This is because of the relatively large negative coefficients that are required to simulate a relative complement operation (revisit Fig \ref{fig2}).

One way to remedy such phenomenon is to consider an asymmetric $\ell_1$ norm defined as

\begin{equation}\label{eq42}
\|\balpha|_{1,w}\triangleq\sum_{\{i:\;\alpha_i\geq 0\}}|\alpha_i|+ \sum_{\{j:\;\alpha_j< 0\}}w|\alpha_j|,
\end{equation}
where $w$ is a constant scalar. When $w\in(0,1)$, this representation allows the appearance of larger negative coefficients by weighting them less in the absolute sum. Additionally, setting $w\gg 1$ promotes positivity on the coefficients which is desirable in some applications that we will consider in next section.

Clearly $\|.|_{1,w}$ is not formally a norm as it violates the positive homogeneity of a norm, however, the corresponding ball is still convex (see Fig \ref{fig6}(b)). A generalization of \cite{van2008probing} reported in \cite{van2011sparse} allows replacing the $\ell_1$ constraint with any convex constraints. As stated in \cite{van2011sparse} for this generalization the infinity norm in (\ref{eq25}) needs to be replaced with its \emph{polar} (reducing to the dual norm when the constraint is a norm), which in this case is the asymmetric infinity norm
\begin{equation}\label{eq43}
\|\balpha|_{\infty,w^{-1}}\triangleq \max_{\begin{subarray}{l}
        i:\;\alpha_i\geq0\\  j:\;\alpha_j<0
      \end{subarray}} \{|\alpha_i|,w^{-1}|\alpha_j|\}.
\end{equation}

\begin{figure}[!htbp]
\centering \includegraphics[width=70mm]{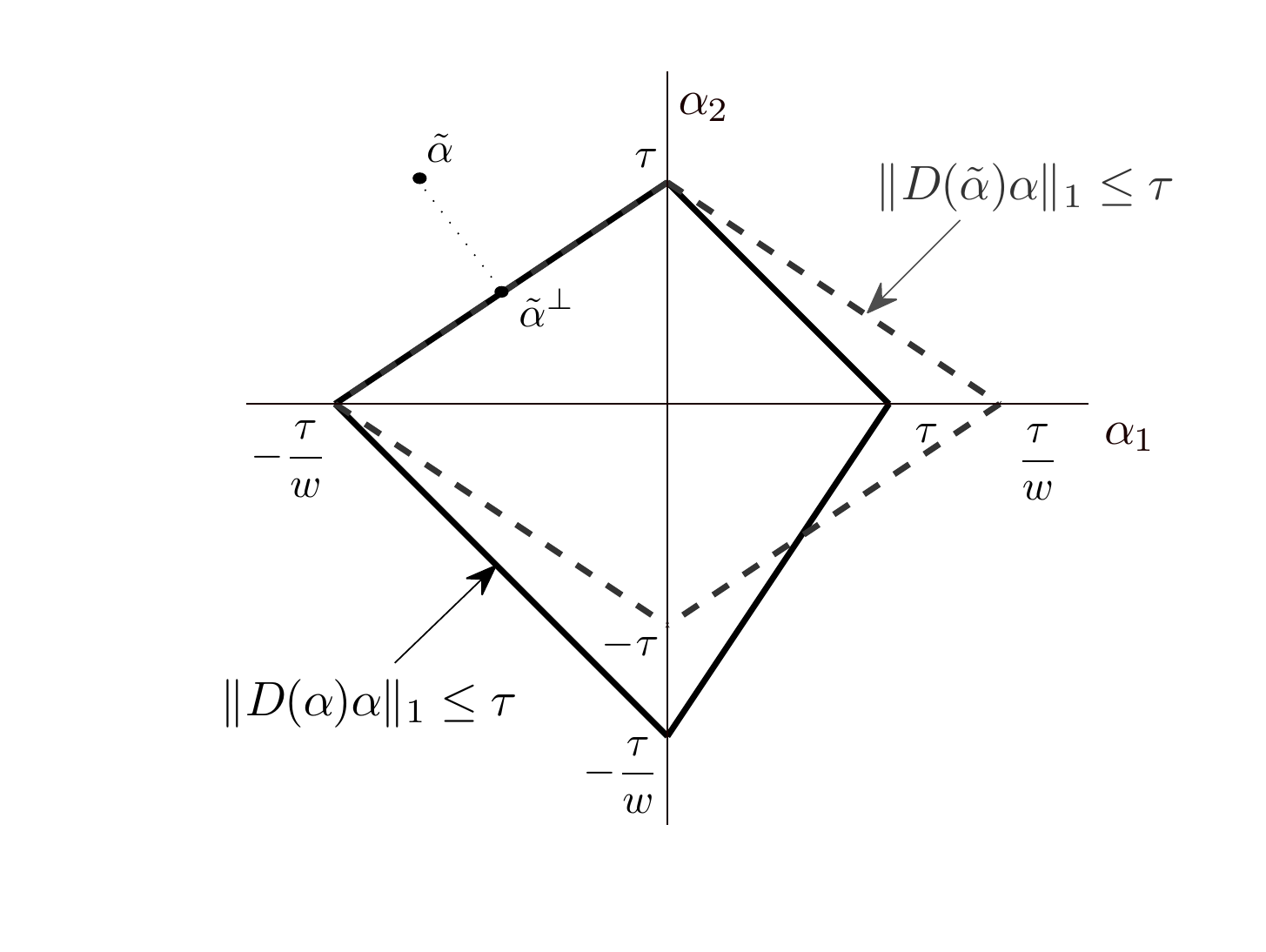}
\caption{Coincidence of the projections onto the asymmetric balls $\|\balpha|_{1,w}\leq \tau$ and $\|\boldsymbol{D}\hspace{-1mm}\mbox{ (}\hspace{-.3mm} \tilde \balpha \hspace{-1mm} \mbox{ )} \balpha\|_1\leq \tau$. The two balls share side in the second quadrant where $\tilde\balpha$ is located. }%
\label{fig6-2}%
\end{figure}

To more conveniently facilitate our method with this feature, we would note that in a lower implementation level, a main component of the SPG method in \cite{van2008probing} is performing iterative projections of the form
\begin{equation}\label{eq44}
\tilde\balpha^\perp= \operatorname*{arg\,min}_{\balpha}\;\;\|\balpha-\tilde\balpha\|_\mathbb{S} \quad s.t. \quad \|\balpha|_{1,w}\leq \tau.
\end{equation}
The asymmetric norm maybe written as $\|\balpha|_{1,w}=\|\boldsymbol{D}\hspace{-.1mm}\mbox{\small (}\hspace{-.3mm} \small\balpha \hspace{-.3mm} \mbox{\small )} \normalsize\balpha\|_1$ where $\boldsymbol{D}\hspace{-.1mm}\mbox{\small (}\hspace{-.3mm} \small\balpha \hspace{-.3mm} \mbox{\small )}$ is a diagonal matrix with entries
\begin{equation}\label{eq44p5}
\boldsymbol{D}\hspace{-.1mm}\mbox{\small (}\hspace{-.3mm} \small\balpha \hspace{-.3mm} \mbox{\small )}_{i,i}= \left\{
     \begin{array}{lr}
       1 &  \alpha_i\geq 0\\
       -w &  \alpha_i<0
     \end{array}.
   \right.
\end{equation}
As we have shown in Appendix B, solution of (\ref{eq44}) coincides with the solution of the weighted $\ell_1$ minimization problem
\begin{equation}\label{eq45}
\min_{\balpha}\;\;\|\balpha-\tilde\balpha\|_\mathbb{S} \quad s.t. \quad \|\boldsymbol{D}\hspace{-.1mm}\mbox{\small (}\hspace{-.3mm} \small\tilde \balpha \hspace{-.3mm} \mbox{\small )} \normalsize\balpha\|_1\leq \tau.
\end{equation}
This fact is well demonstrated in Fig \ref{fig6-2}. Based on this argument we may replace the projection onto an asymmetric ball problem with a weighted $\ell_1$ minimization, which is straightforward.

\section{Examples}\label{exsec}
In this section we examine the proposed technique in various imaging applications. For the examples presented, we use the approximate Heaviside function suggested in \cite{aghasi2011parametric} as
\begin{equation}
H_\epsilon(x)= \left\{
     \begin{array}{lr}
       1 &  x> \epsilon\\
       0 &  x<\epsilon\\
       \frac{1}{2}+\frac{x}{2\epsilon}+\frac{1}{2\pi}\sin(\frac{\pi x}{\epsilon})& |x|\leq \epsilon,
     \end{array}
   \right.
\end{equation}
which gives rise to a compactly supported approximation of $\delta_\epsilon(x)=H_\epsilon'(x)$. As discussed in \cite{aghasi2011parametric} in the context of inverse problems, this choice results in a form of narrow-banding in the evolution of the parameters, that will also be discussed briefly in segmentation examples. The parameter $\epsilon$ determines the width of the narrow-band region, whereas smaller values of $\epsilon$ brings sharper transitions into the reconstructions at the expense of slower evolution \cite{aghasi2011parametric}. The choice of the lifting parameter $c$ is quite arbitrary as the proposed shape-based problem is relative, i.e., simultaneously scaling $c$ and $\alpha_i$ coefficients in (\ref{eq17}) does not change the zero level set shape. In other words, taking larger values of $c$ would tend to larger values of $\alpha_i$ in the reconstructions. Since in general the shapes in the dictionary may significantly vary in size, to have knoll basis terms of comparable magnitudes, we normalize each knoll to its maximum value. For the majority of examples we simply take $c=0.1$ and $\epsilon=0.05$ and initialize the algorithm from a rather random state. Based on the scaling property between $c$ and final $\alpha_i$ values, $\tau_{mx}$ may be taken as a large multiple of $c$, e.g., $\tau_{mx}=50c$.
\subsection{Image Segmentation}
As discussed in Section \ref{sec-Im2}, for a binary Chan-Vese segmentation the energy functional finds the form of (\ref{eq34}) which maybe written as the sum of internal and external energies:
\begin{equation}\label{eq46}
\mathcal{E}(\balpha)=\|\mathcal{G}_{in}(\balpha)\|_\mathbb{S}^2+ \|\mathcal{G}_{ex}(\balpha)\|_\mathbb{S}^2.
\end{equation}
Accordingly a reasonable representation for $\mathcal{G}(\balpha)$ would be
\begin{equation}\label{eq47}
\mathcal{G}(\balpha)= \mathcal{G}_{in}(\balpha)+ \ii\mathcal{G}_{ex}(\balpha),
\end{equation}
where $\ii=\sqrt{-1}$. The underlying Hilbert space would be $\mathbb{S}=L^2(D)$ with the corresponding inner product
\begin{equation}\label{eq48}
\langle s_1(x), s_2(x) \rangle_{L^2(D)}=\int_D s_1(x)\overline{s_2(x)}\mbox{d}x,
\end{equation}
where $\overline{s_2(x)}$ represents the complex conjugate of $s_2(x)$. Clearly, applying our proposed algorithm to the segmentation problem requires having $\mathcal{G}'(\balpha)[.]$ and its adjoint form. For a given vector $\boldsymbol{\eta}\in\mathbb{R}^{n_d}$, a formal derivation of (\ref{eq47}) with respect to $\balpha$ and applying the resulting Jacobian operator to $\boldsymbol{\eta}$ yields
\begin{equation}\label{eq49}
\mathcal{G}'(\balpha)\boldsymbol{\eta} = f(x,\balpha)\sum_i^{n_d}\eta_i \psi_{\mathcal{S}_i}(x),
\end{equation}
where
\begin{equation}\label{eq50}
f(x,\balpha)= \left\{
     \begin{array}{lr}
       \frac{1}{2}\Big(\sqrt{\frac{r_{in}(u)}{H_\epsilon(\phi)}}-\ii \sqrt{\frac{r_{ex}(u)}{1-H_\epsilon(\phi)}}\Big)\delta_\epsilon (\phi) &  |\phi|<\epsilon\\
       0 &  else.
     \end{array}
   \right.
\end{equation}
From a numerical perspective, the linear operator $\mathcal{G}'(\balpha)$ is a matrix with columns $f(x,\balpha)\psi_{\mathcal{S}_i}(x)$. Based on this argument, the adjoint operator may be easily deduced. For a function $\psi\in L^2(D)$ the linear operation $\mathcal{G}'^{*}(\balpha)\psi$ produces an array of size $n_d$ where
\begin{equation}\label{eq51}
[\mathcal{G}'^{*}(\balpha)\psi]_i= \Big\langle \overline{f(x,\balpha)}\psi_{\mathcal{S}_i}(x), \psi(x)\Big\rangle_{L^2(D)}.
\end{equation}
To maintain generality, we preferred providing the implicit forms (\ref{eq49}) and (\ref{eq51}) for $\mathcal{G}'(\balpha)$ and $\mathcal{G}'^{*}(\balpha)$, however, for pixelated images of moderate size a matrix may be assigned to each linear operator.

Equation (\ref{eq50}) shows that $f(x,\balpha)$ is a compactly supported function. Similarly, $\psi_{\mathcal{S}_i}(x)$ is compactly supported, and therefore when $\mbox{supp}\big(\delta_\epsilon(\phi)\big)\cap \mbox{supp}(\psi_{\mathcal{S}_i})=\emptyset$, the corresponding columns of $\mathcal{G}'(\balpha)$ would become zero. In determining a descent direction, such columns may be neglected as they would not affect the result. In other words, in iteratively updating the knoll coefficients $\alpha_i$, the only knolls that are updated are the ones that intersect the narrow-band $\delta(\phi)$ at that iteration. This numerical advantage is basically the \emph{narrow-banding} feature associated with the proposed parametric representation and extensively described in \cite{aghasi2011parametric} in the context of inverse problems.

The technical details presented may be employed to address variants of Chan-Vese segmentation in different applications. Here we consider two applications that could be addressed through a sparse shape recovery.
\subsubsection{Image Segmentation with Missing Pixels}
A challenging problem associated with image segmentation arises when some of the pixels are missing or part of the image is occluded. This limits the analysis to only a subset of image pixels, while the segmentation needs to be performed globally. In this example we show that prior information about the geometry of objects in the image can lead to perform a completion to missing pixels when a sparse shape recovery is considered.

\begin{figure}[!htbp]
\hspace{-3mm}
\includegraphics[width=139.5mm]{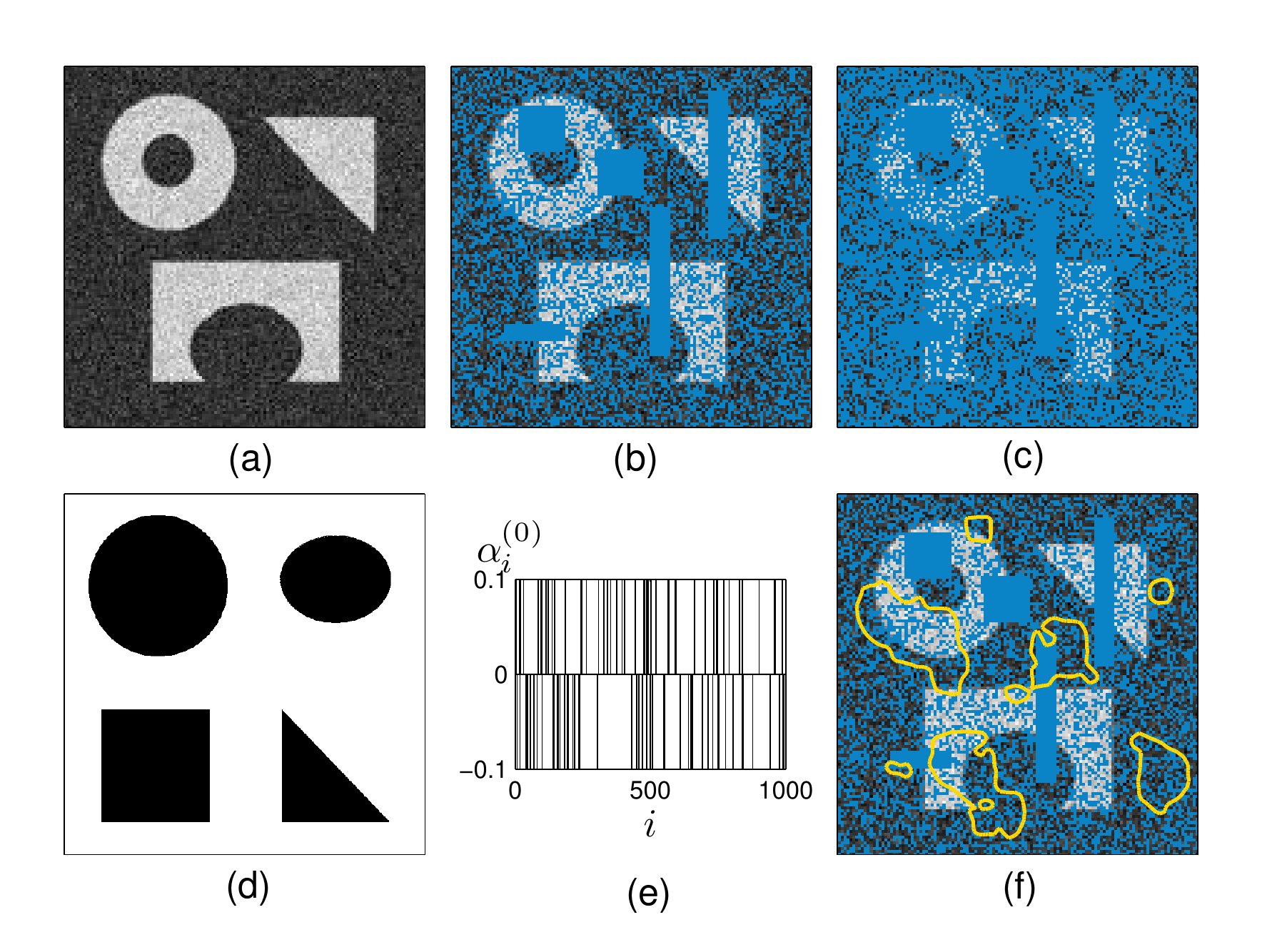}
\caption{(a) Reference image (b) A test image with more than 50\% of the pixels missing (blue regions show the missing regions) (c) A test image with more than 80\% of the pixels missing (d) Shapes used to build up the dictionary: 1000 instances of these four shapes with different sizes are placed throughout the imaging domain (e) Initial values of $\alpha_i$ for $i=1,\cdots 1000$ (f) The resulting initial contour for the given initial $\balpha$. The segmentation results are shown in Fig \ref{fig8}.}\label{fig7}
\end{figure}

\begin{figure*}[!htbp]
\hspace{-5mm}
\begin{tabular}{ccc}
\includegraphics[width=41mm]{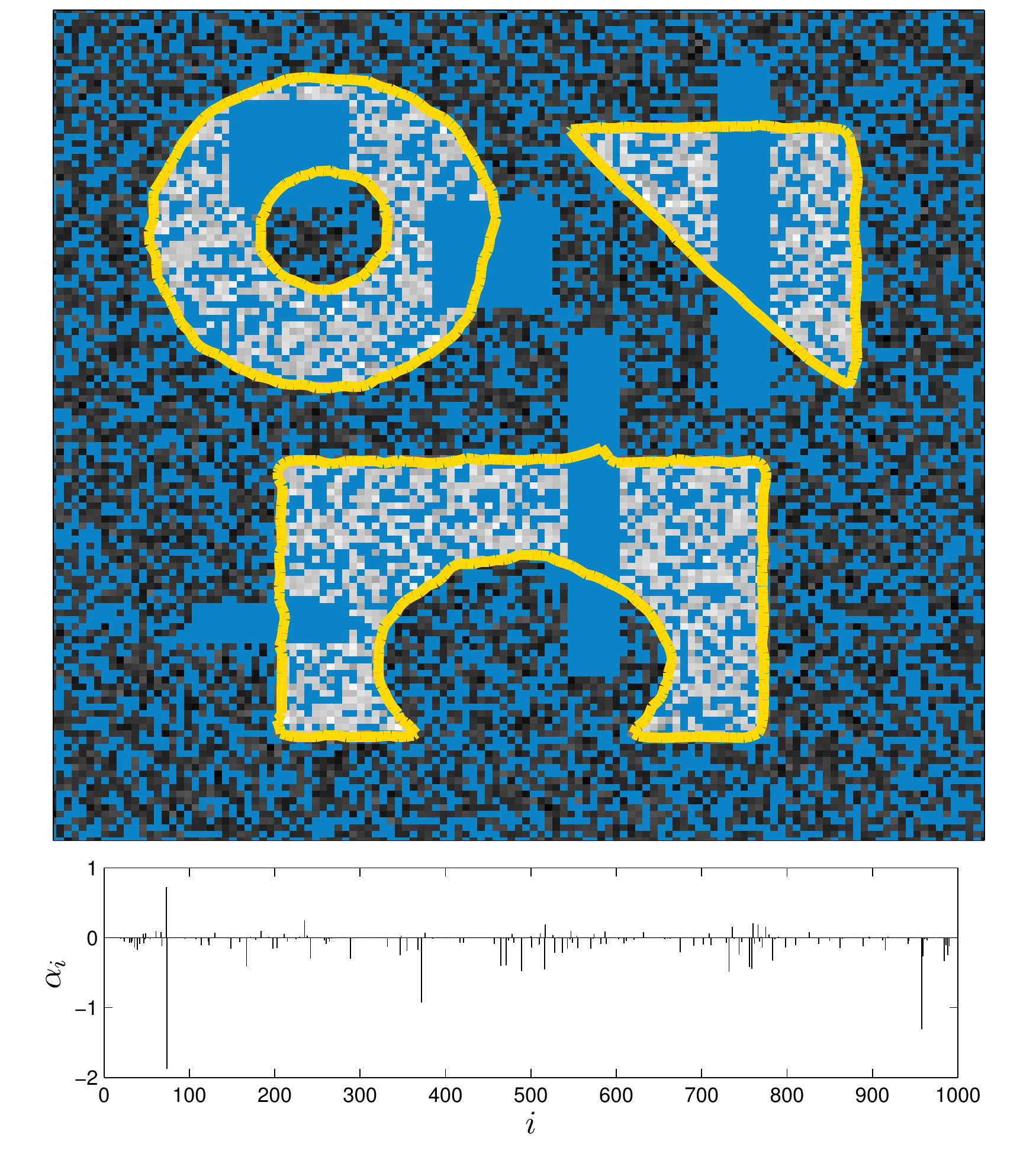} &
\includegraphics[width=41mm]{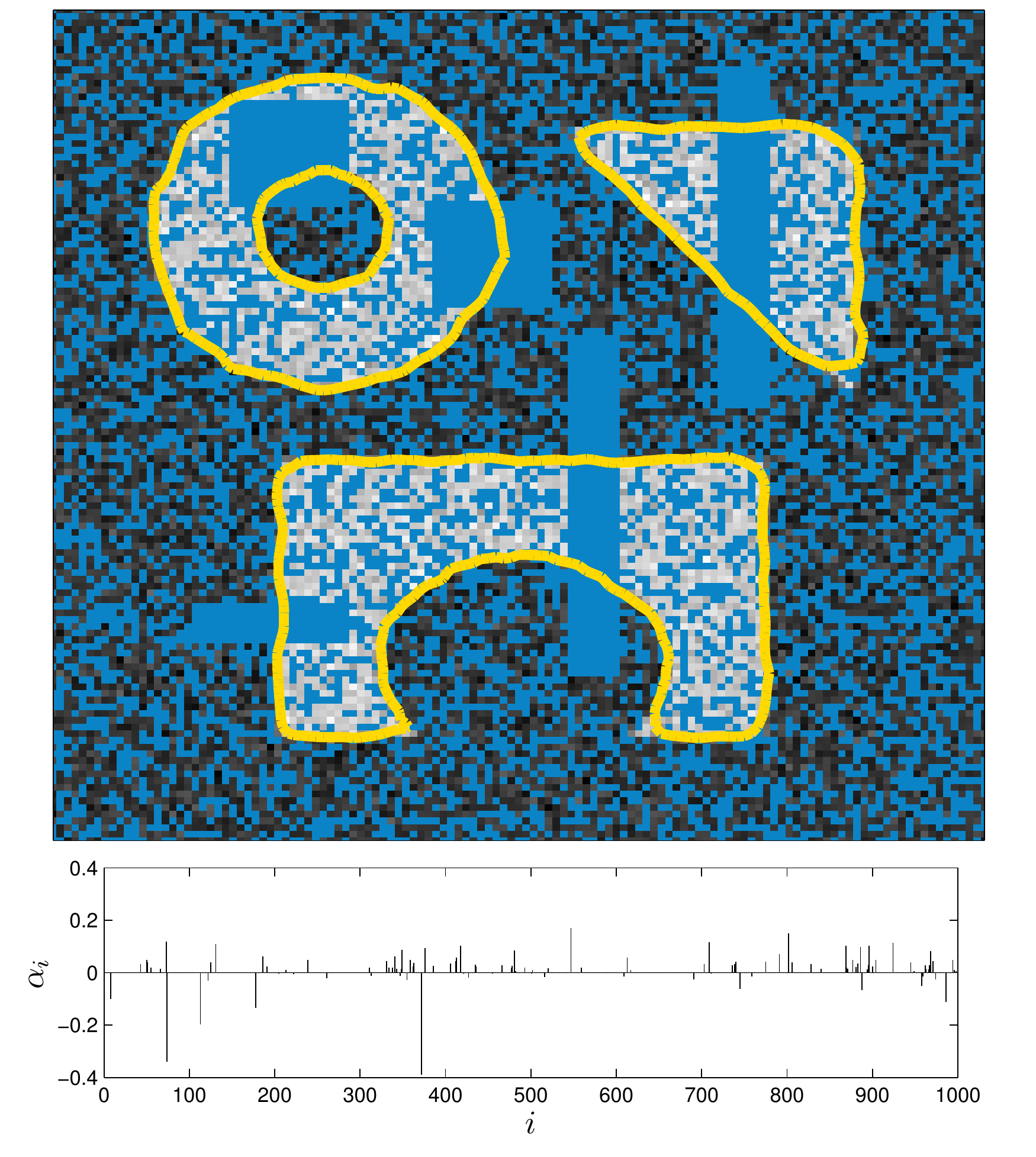}&
\includegraphics[width=41mm]{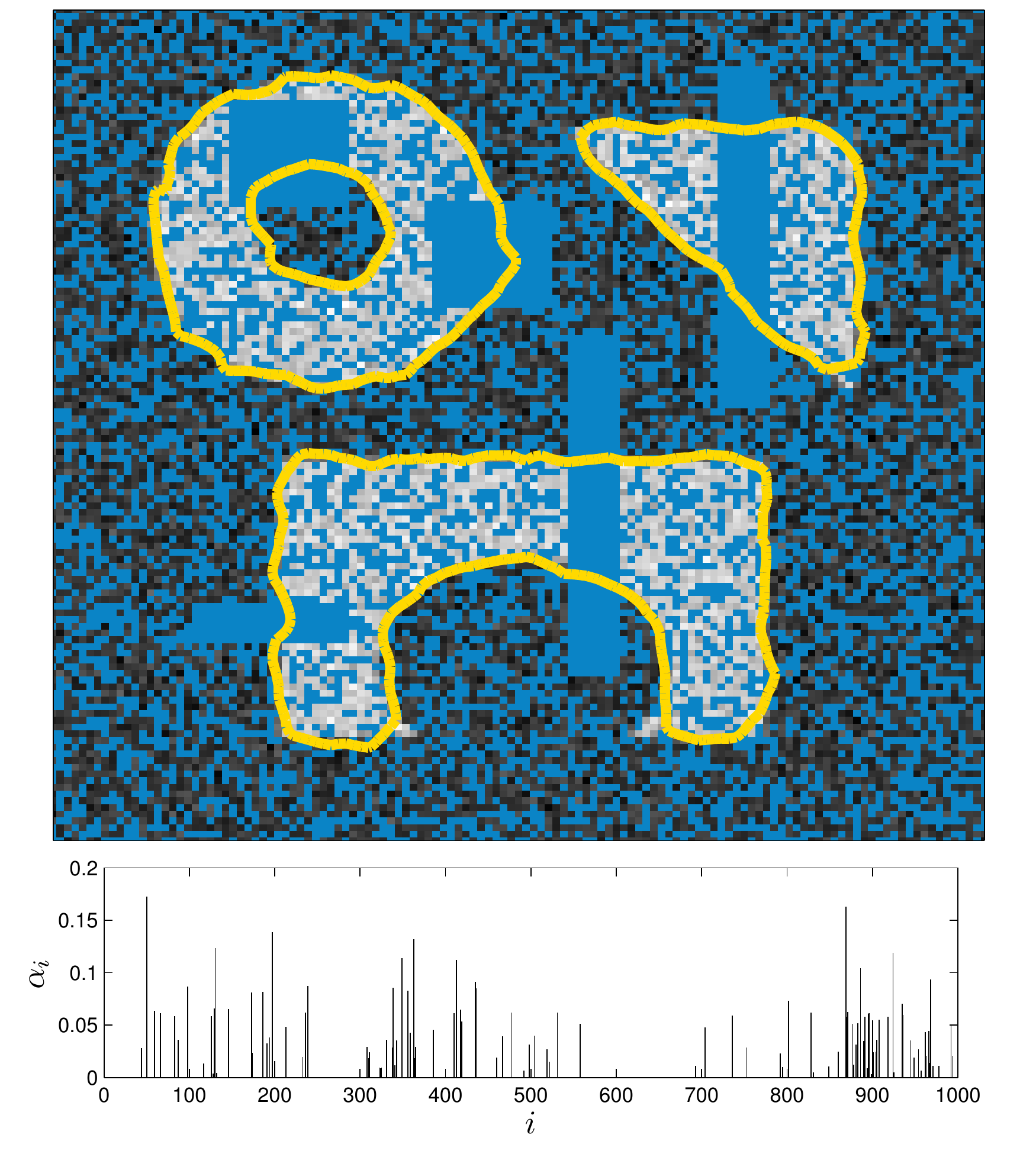} \\[-.2cm]
\includegraphics[width=41mm]{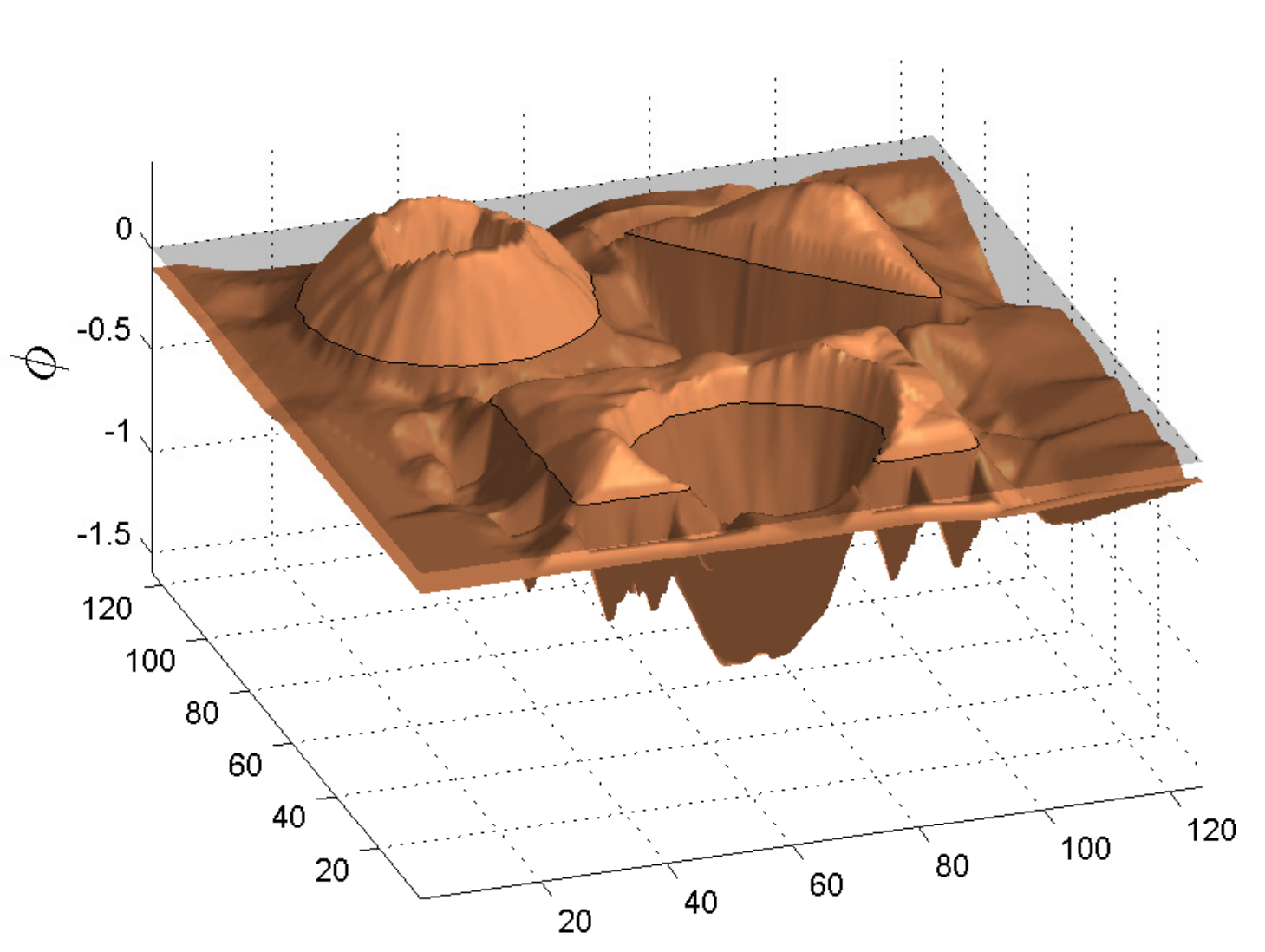} &
\includegraphics[width=41mm]{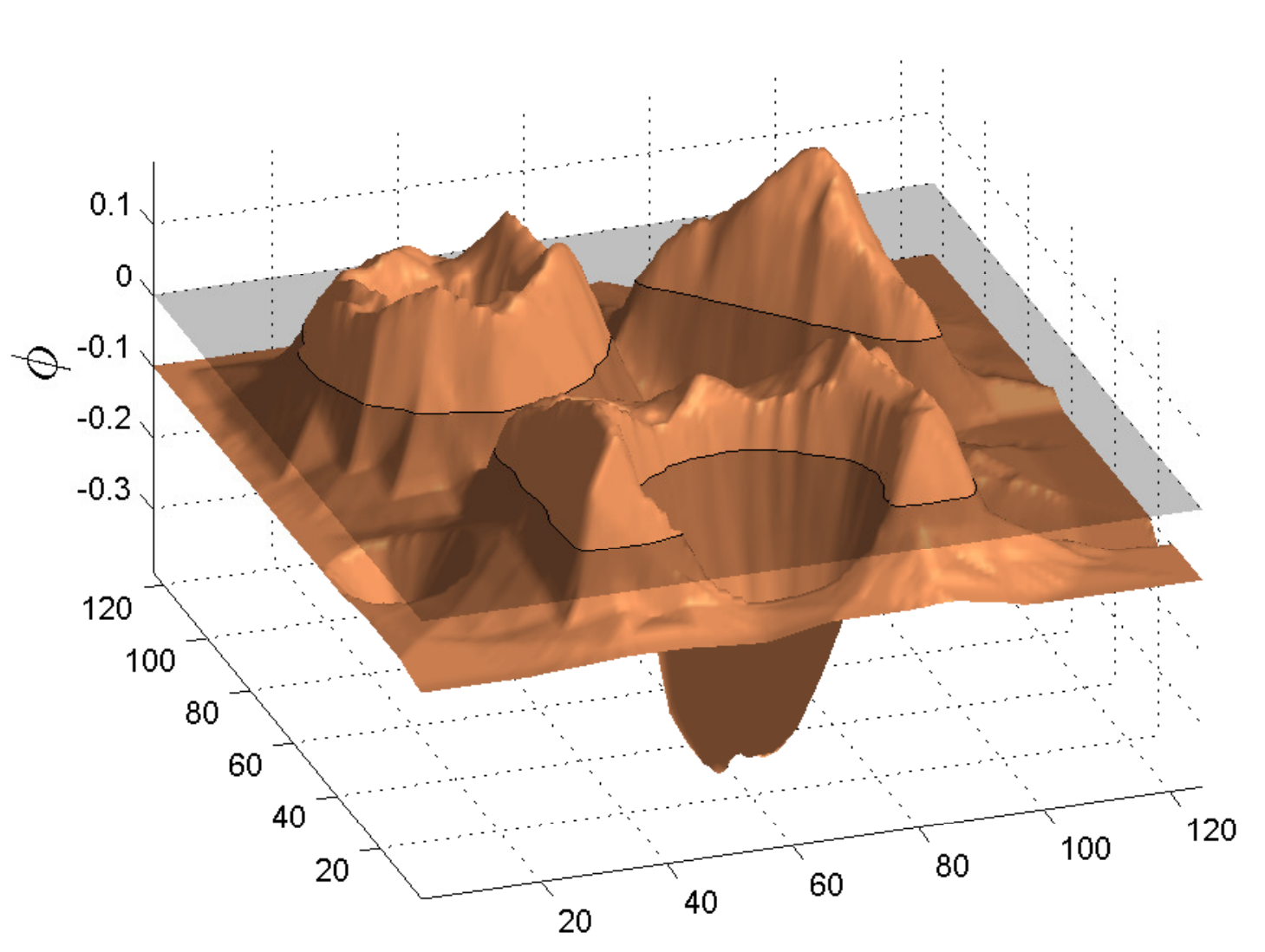}&
\includegraphics[width=41mm]{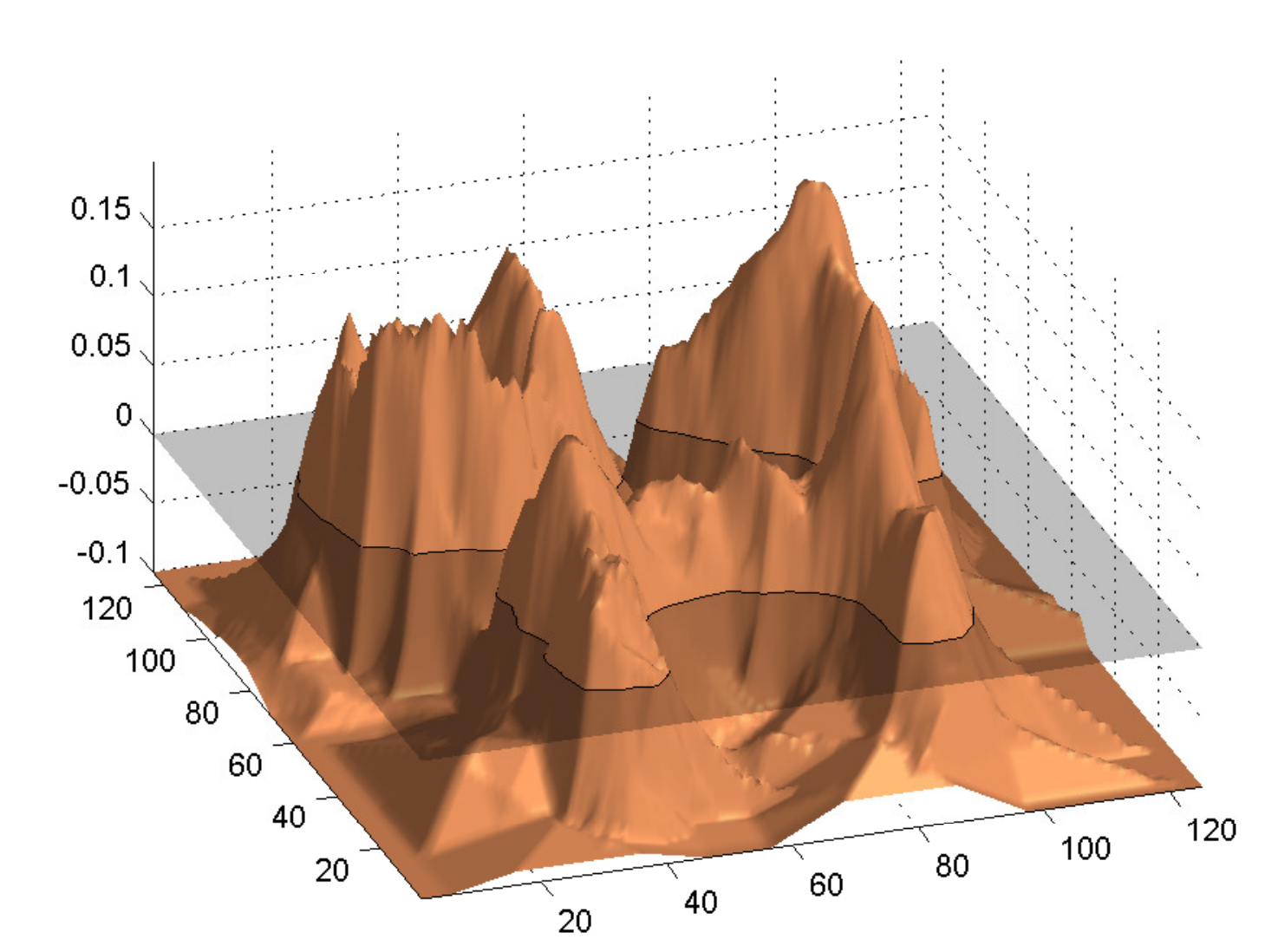} \\
\includegraphics[width=41mm]{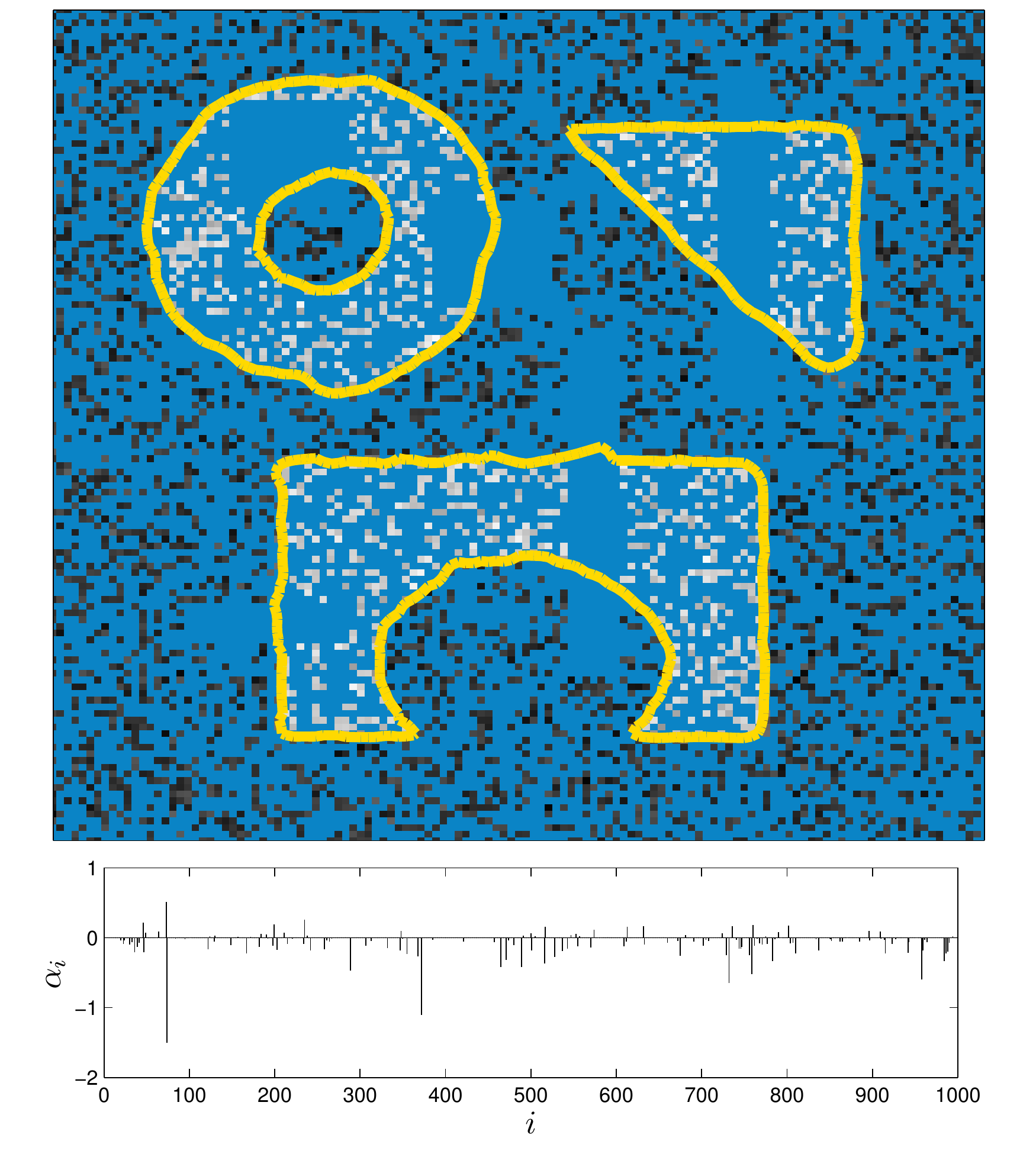} &
\includegraphics[width=41mm]{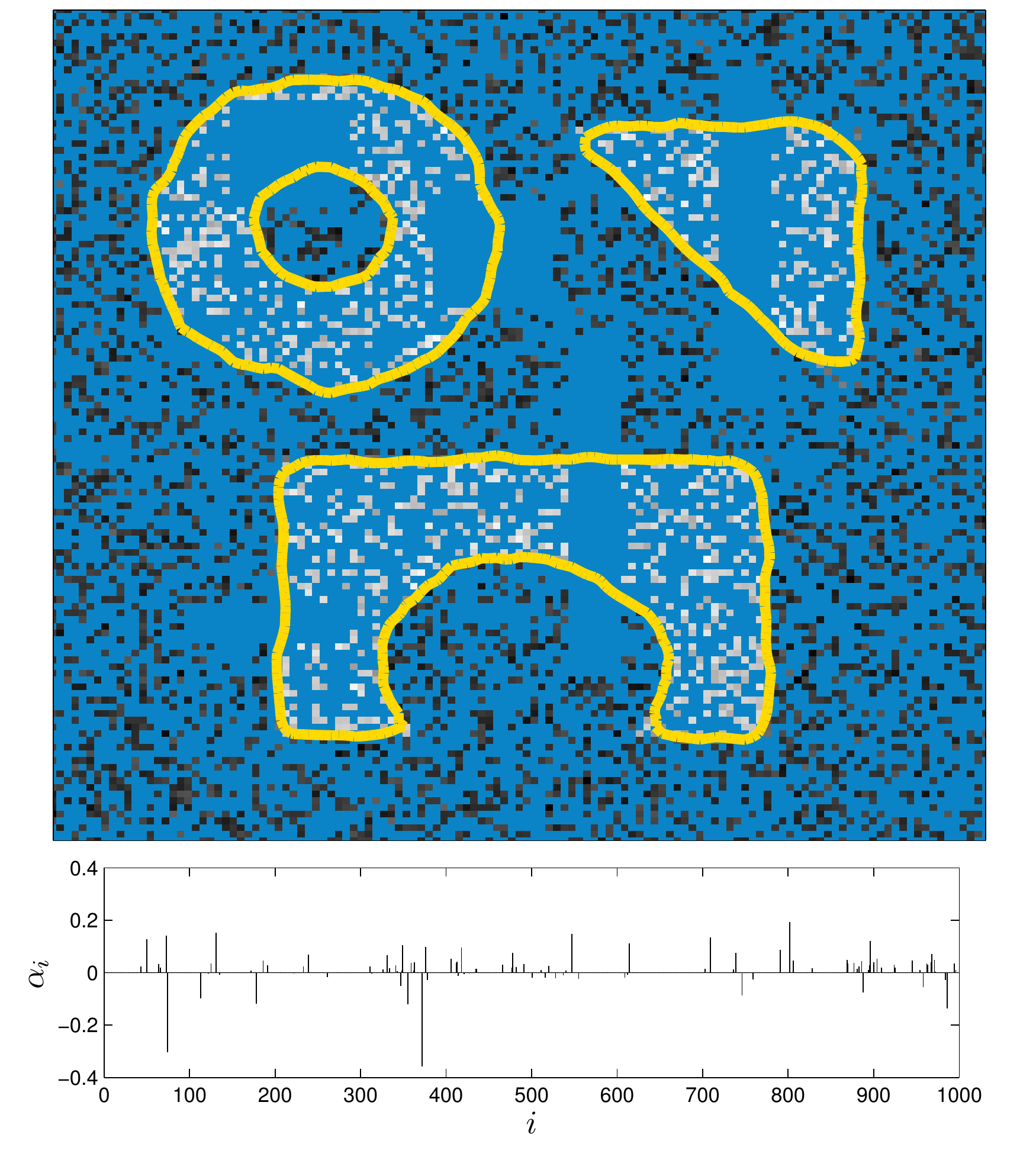}&
\includegraphics[width=41mm]{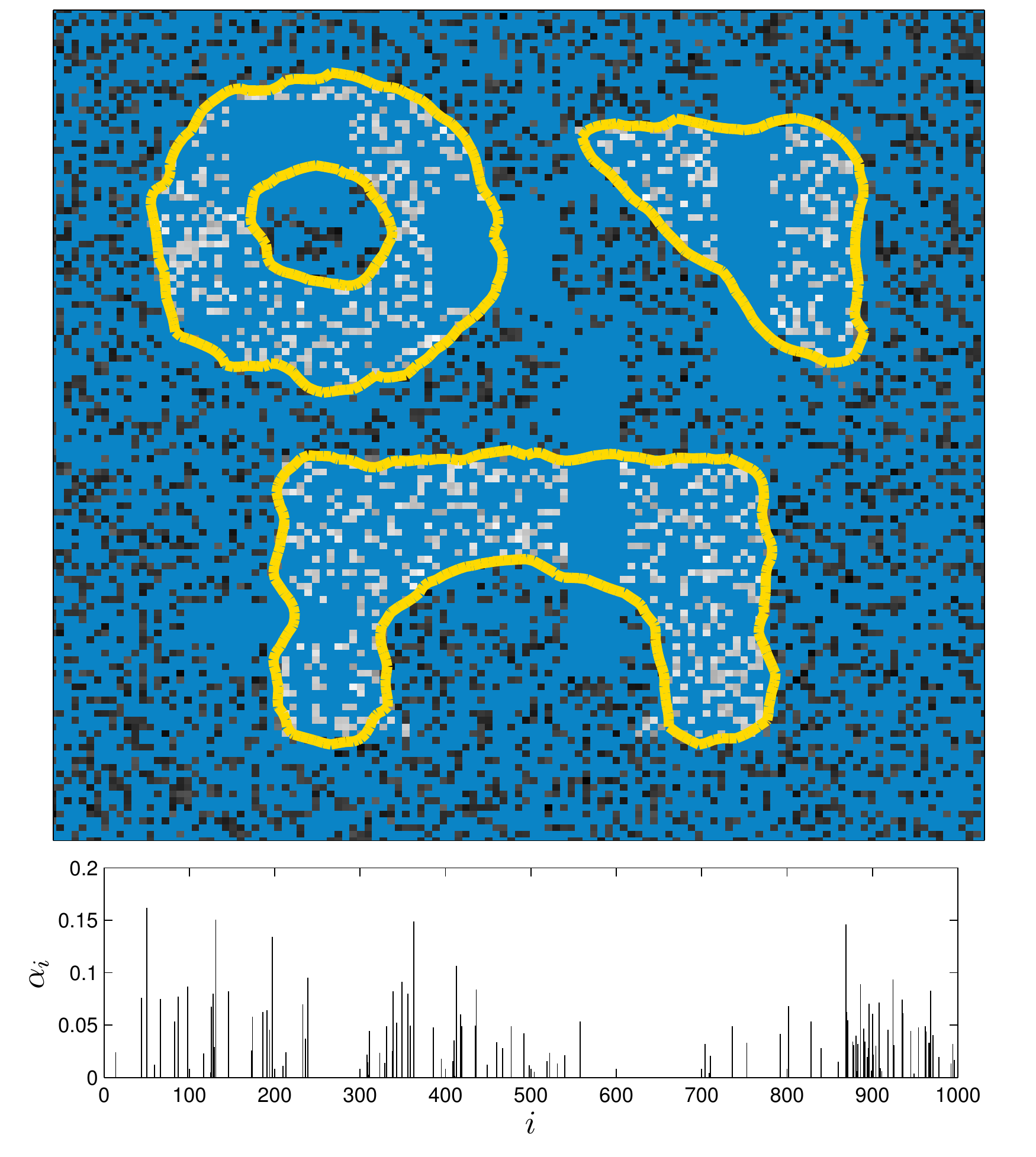}
\end{tabular}
\caption{First row shows the results of segmentation on 50\% missing-pixel image, from left to right for $w=0.1$, $w=1$ and $w=10$. Below every result the final sparse vector $\balpha$ is shown. Second row are the level set functions corresponding to each segmentation above it. Third row shows an identical scenario as the first row applied to the image with 80\% missing-pixel}\label{fig8}
\end{figure*}

Fig \ref{fig7} shows a noisy reference image that is occluded with rectangular patches. Additionally, random pixels of the image are discarded resulting in overall pixel loss of 50\% and 80\% shown in Figures \ref{fig7}(b) and \ref{fig7}(c). These two figures will be the subject of our segmentation.

Our prior information about the geometry of objects in the image is reflected in the choice of dictionary elements. To build up the dictionary we make use of four basic shapes: a circle, square, triangle and an ellipse, shown in Fig \ref{fig7}(d). We would note that the triangle in the reference image is the upper diagonal portion of a square while the one used in the dictionary is a lower portion.

A slight modification of the formulation presented above needs to be considered for this limited observation problem. More specifically, considering $D'\subset D$ to be the available portion of the imaging domain $D$, the Hilbert space considered would be $L^2(D')$ where the observation takes place. The knoll functions, however, are constructed and defined globally in $D$.

By considering different size and placements of the four basic shapes, a total of $n_d=1000$ shapes are assigned to the dictionary. To initialize the algorithm we randomly apply weights of $+1$ and $-1$ to 100 knolls as shown in Fig \ref{fig7}(e). This assignment results the initial shape shown in Fig \ref{fig7}(f). A basic Chan-Vese segmentation cost is considered for which the values $\tilde u_{in}$ and $\tilde u_{ext}$ are occasionally updated following the formulation in \cite{chan2001active}.

The first column of Fig \ref{fig8} shows the successful segmentation results on the two test images that are obtained for $w=0.1$. Although significant parts of the images are missing, in both cases the proposed algorithm has made a reasonable segmentation job. To highlight the advantage of considering an asymmetric norm, results are also shown for the basic $\ell_1$ case of $w=1$ and an inappropriate choice of $w=10$ which basically abandons set minus operation. It can be seen that although considering smaller values for $w$ increases the size of feasible region and requires solving the problem in a larger domain, it pays off by promoting the set minus operation and using a sparser set of shapes in a more efficient manner. This contrast is clearly observable by comparing the results of first and third columns. Basically, in the third column knolls are pushed to take positive weights which degrades the reconstruction by exclusively promoting the union operation. Thanks to the appearance of the set minus operation, the segmentation results in the first column have sharper corners and smoother sides closer to the truth while using relatively less number of shapes in the dictionary. It is worth noting that all segmentations converged the steady state in less than 15 iterations.

\subsubsection{Text Recognition: Breaking a Basic Captcha}
Captcha (Completely Automated Public Turing test to tell Computers and Humans Apart) is a well-known test in the world of computers that ensures that response to a query is generated by a human \cite{von2003captcha}. For this test, the client is asked to read and type the word in an image, where the characters are placed in an unusual manner not recognizable for the machine. Since the underlying components of the image are still characters, this maybe considered as a challenging sparse shape reconstruction problem. The underlying shapes in the dictionary may constitute a dense set of possibilities for the character shapes that may appear in a Captcha image. Of course a denser dictionary increases the chances of identifying the characters correctly.

\begin{figure}[!htbp]
\centering
\includegraphics[width=126mm]{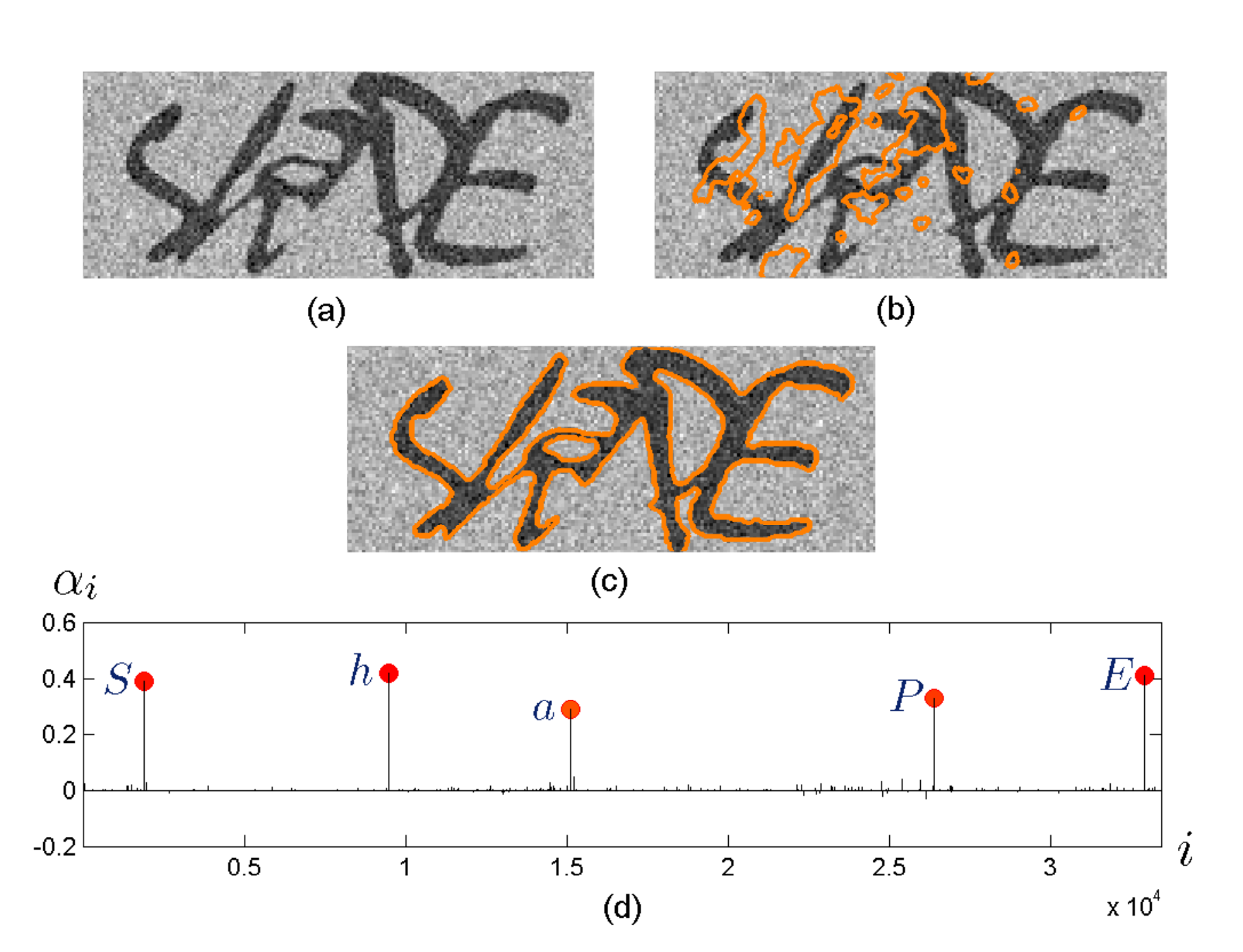}\vspace{-2 mm}
\caption{(a) A Captcha image (b) The initial contour corresponding to the random initialization of the $\alpha_i$ coefficients (c) The segmentation result shown with colored contour (d) The values of reconstructed weights and an indication of the letter each index corresponds to }\label{fig9}
\end{figure}

As a proof of concept, in this example we consider an attack on the rather basic Captcha image of size $63\times 160$ pixels shown in Fig \ref{fig9}(a) which represents the word ``\emph{ShaPE}''. The letters in the image may take different case, size, rotations and overlap as the case for our Captcha image. To build up the character dictionary we assume knowing the font type a priori and consider all 52 uppercase and lowercase letters of the English alphabet. The displacement possibilities considered for each character are 64 points on a regular grid in the image. At every point, five different rotations of the character and two different font sizes are considered. This setup leads to $n_d=52\times 64 \times 5\times 2=33280$ shapes in the dictionary. Of course for a more realistic scenario, font possibilities, deformations and more number of size, displacement and rotational variations may be considered for every letter. We however keep the problem small to be tractable with a desktop computer. Since the main purpose of this problem is identifying shapes in the dictionary that appear somewhere in the image, $w$ may be chosen to be large (10 for our simulations) to push the algorithm on only considering the union operation.

As the first experiment we consider the case that exact character shapes in the Captcha are present in the dictionary. As before, the algorithm is initialized with positive and negative weights for a random subset of $\alpha_i$ coefficients which corresponds to the initial shape contour in Fig \ref{fig9}(b). The classic Chan-Vese cost functional is again considered with the texture parameters updated occasionally. The segmentation result after 12 iterations is shown in Fig \ref{fig9}(c) and in Fig \ref{fig9}(d) we have shown the coefficients values. A simple index map relating the weight indices to the alphabetical representation of the shape indicates that the reconstructed image is composed of five main characters, precisely matching with the word inside the image.

\begin{figure}[!htbp]
\centering
\includegraphics[width=126mm]{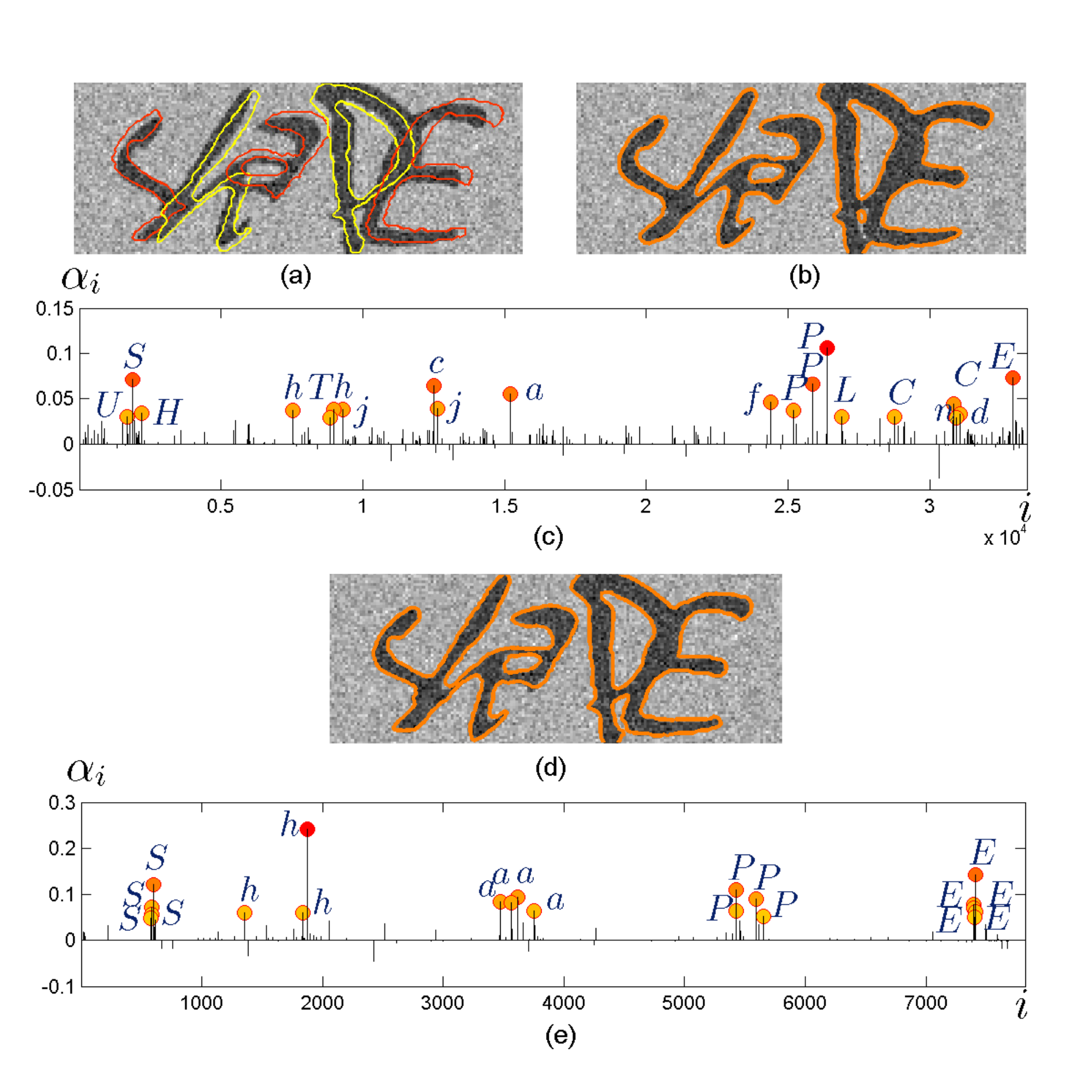}\vspace{-5 mm}
\caption{(a) A Captcha image where the letters in the image do not match with the elements of the dictionary; the closest elements of the dictionary are shown with colored contours (b) A first stage reconstruction results (c) Reconstructed weights corresponding to the first stage segmentation; the top 20 candidates are indicated with the corresponding letter (d) Segmentation results after refining the dictionary according to the first stage results (e) The weights corresponding to the second stage reconstruction indicating the letter each index corresponds to }\label{fig10}
\end{figure}
For the second example we consider the Captcha in Fig \ref{fig10}(a) where the letters in the image do not precisely match the dictionary elements. The closest elements of the dictionary are depicted with colored contours on the same image. This problem may be still solved using a multi-stage refinement strategy. More specifically, a first segmentation attempt shown in Fig \ref{fig10}(b) with the underlying weights shown in Fig \ref{fig10}(c) does not provide a sparse solution with prominent weights as before. Instead, in this case there are several candidates, the top 20 of which are shown in Fig \ref{fig10}(c). Based on the first stage results, a new dictionary maybe built up, with much a much smaller number of elements, only containing denser size, rotation, displacement variants of the shapes listed in the top 20 list. Performing the segmentation process again, this time provides us with the weights shown in Fig \ref{fig10}(e) where the striking weights only correspond to the letters in the Captcha. The ``\emph{echo}'' effect corresponding to appearance of multiple prominent weights of the same letter is due to the close distance between the shapes in the new dictionary. Of course when the indexing mechanism takes into account the character positions, all such echoes correspond to close positions and maybe unified as a single final conclusion about the letter in that region of the image.

\subsection{Medical Imaging: X-ray Computed Tomography (CT)}\label{insec1}
As a well known linear inverse problem, in this section we examine the method in reconstructing CT images. For a mono-energetic CT, X-ray photons are transmitted through the test medium and measured at the opposite side. If the medium has an attenuation profile $\mu(x)$ the number of photons measured would ideally be
\begin{equation}\label{eq52}
\lambda_m= \lambda_T \exp\big(-\int_{\mathcal{L}_m}\!\!\!\mu(x)\mbox{d}x\big), \quad m=1,\cdots,M,
\end{equation}
where $\lambda_m$ is the photon count measured at $m$-th receiver, $\lambda_T$ is the blank scan photon count and $\mathcal{L}_m$ is the line through which the ray travels. An easier way of interpreting the measurements is reading the values
\begin{equation}\label{eq53}
 v_m= -\log \frac{ \lambda_m}{\lambda_T}= \int_{\mathcal{L}_m}\!\!\!\mu(x)\mbox{d}x,
\end{equation}
at each sensor which basically relates the attenuation map to the measurements via a Radon transform. In practice however, the measurements are corrupted with Poisson noise, i.e.,
\begin{equation}\label{eq54}
\tilde \lambda_m=\mbox{Pois}(\lambda_m)
\end{equation}
and the true measurements are $\tilde v_m=-\log(\tilde \lambda_m/\lambda_T)$. A second order approximation to the log likelihood function indicates that \cite{sauer1993local}
\begin{equation}\label{eq55}
\log p(\tilde{\boldsymbol{v}}|\mu)\approx -\frac{1}{2}(\tilde{\boldsymbol{v}}-\mathcal{R}\mu)^T \boldsymbol{D} (\tilde{\boldsymbol{v}}-\mathcal{R}\mu) + h(\tilde{\boldsymbol{v}})
\end{equation}
where $\boldsymbol{D}=\mbox{diag}(\tilde \lambda_1, \cdots, \tilde \lambda_M)$, $\mathcal{R}$ is the linear Radon type transform that maps $\mu$ to the ideal measurements and $h(.)$ is a function dependent only on the data vector $\tilde{\boldsymbol{v}}$. To apply the proposed algorithm to this modality the attenuation profile is parameterized as $\mu(x,\balpha)$, and the residual operator is written as $\mathcal{G}(\balpha)=\mathcal{R}\mu(x,\balpha)- \tilde{\boldsymbol{v}}$. For a maximum likelihood estimate of the parameters, the underlying cost takes the form $\|\mathcal{G}(\balpha)\|_{\boldsymbol{D}}^2$ where the inner product in the measurement space is defined as $\langle \boldsymbol{s}_1, \boldsymbol{s}_2 \rangle_{\boldsymbol{D}}=\boldsymbol{s}_1^T \boldsymbol{D} \boldsymbol{s}_2$ for $\boldsymbol{s}_1, \boldsymbol{s}_2 \in \mathbb{R}^M$. Based on the linearity of the CT model the Jacobian operator applied to a vector $\boldsymbol{\eta}$ may be written as
\begin{align}
\mathcal{G}'(\balpha)\boldsymbol{\eta}=\sum_i{\eta_i \frac{\partial \mathcal{G}(\balpha)}{\partial \alpha_i}}=\sum_i \eta_i \mathcal{R}\frac{\partial \mu}{\partial \alpha_i}=\mathcal{R}\sum_i \eta_i \frac{\partial \mu}{\partial \alpha_i}\label{eq56}
\end{align}
and the adjoint operator applied to a vector $\boldsymbol{\psi}\in \mathbb{R}^M$ takes the form of
\begin{equation}\label{eq57}
\mathcal{G}'^*(\balpha)\boldsymbol{\psi}= \langle \mathcal{R}\frac{\partial \mu}{\partial \alpha_i}, \boldsymbol{\psi}\rangle_ {\boldsymbol{D}}.
\end{equation}
Equations (\ref{eq56}) and (\ref{eq57}) are in general the key components of applying Algorithm 1 to this problem.
\begin{figure}[!htbp]
\centering
\includegraphics[width=135mm]{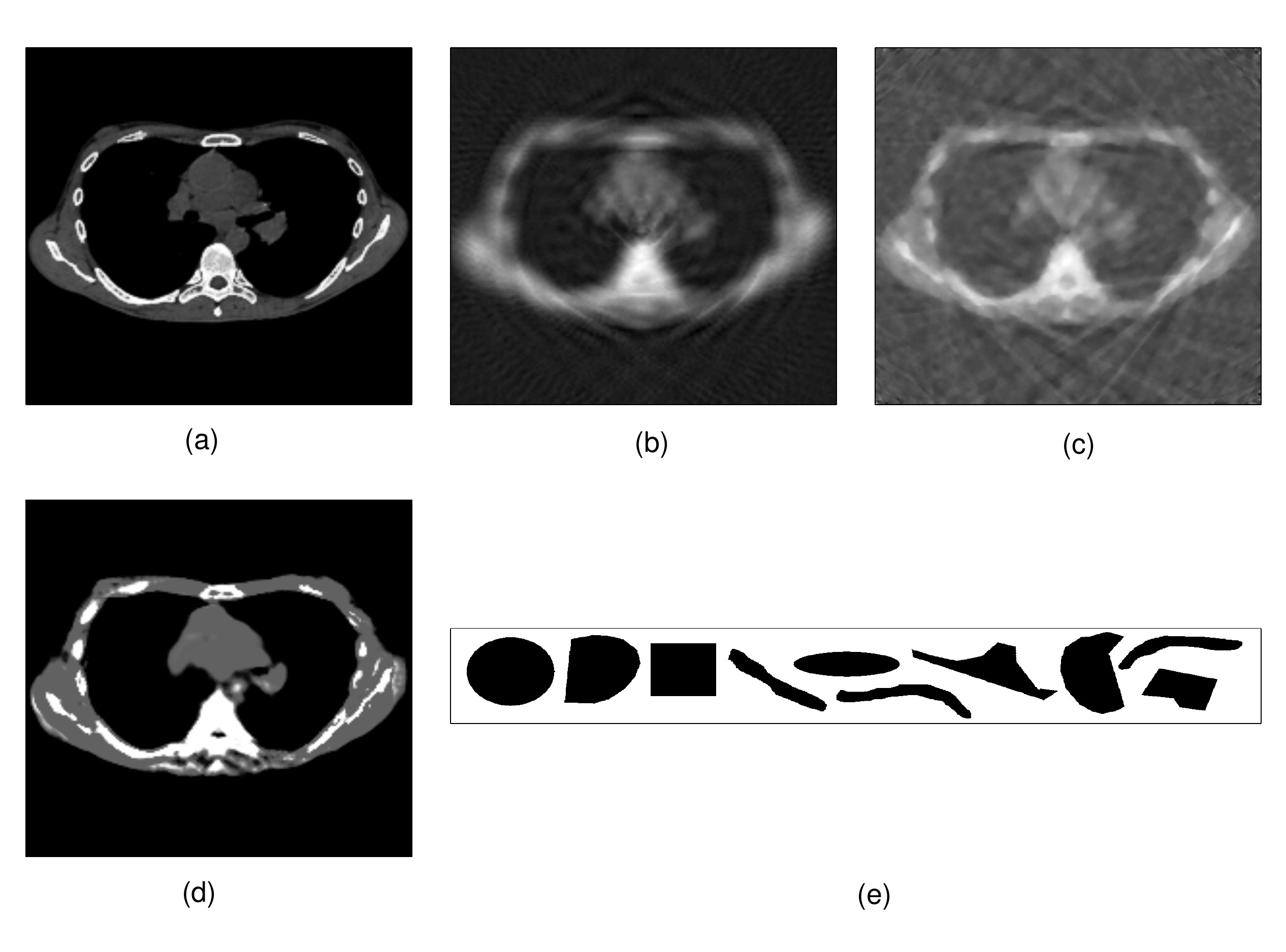}\hspace{-10mm}
\caption{(a) The true attenuation profile of a chest image (b) Inversion using FBP (c) TV inversion (d) Sparse shape inversion results (e) The basic shapes used in the dictionaries}\label{fig11}
\end{figure}

Fig \ref{fig11}(a) shows the true attenuation profile of a $200\times 200$ pixel chest test image. To generate the X-ray data the rays are emitted at 60 equispaced angles between 0 and 180 degrees. At each angle an average of 38 photon rays are transmitted through the image. A total of 2280 measurements are obtained among which 1499 are the ones crossing the chest and containing useful information. The blank scan photon count is $\lambda_T=4\times 10^6$ and the average photon energies are 50 KeV. Beside considering the Poisson noise, 1\% zero mean Gaussian noise is also added to the data to model the inaccuracy in the measurement readings. This problem setup poses a challenging ill-posed problem which is rather hard to approach. In Fig \ref{fig11}(b) we have shown the filtered back projection (FBP) results for this data set. Fig \ref{fig11}(c) also shows the reconstruction results using a total variation (TV) approach employing the $\ell_1$-magic package \cite{candes2005l1}.

For the purpose of shape representation we consider a multi-phase level set approach \cite{vese2002multiphase}. More specifically to invert for both the soft tissue and bone geometries along with the blank space we consider using two level sets $\phi_1(x,\balpha^{(1)})$ and $\phi_2(x,\balpha^{(2)})$ as
\begin{align*}
\mu(x,\balpha^{(1)},\balpha^{(2)})&=\mu_a+(\mu_s-\mu_a)H_\epsilon \big (\phi_1(\balpha^{(1)})\big)\\&+ (\mu_b-\mu_s)H_\epsilon \big (\phi_1(\balpha^{(1)})\big)H_\epsilon \big(\phi_2(\balpha^{(2)})\big)
\end{align*}
where $\mu_a$, $\mu_s$ and $\mu_b$ are respectively the average attenuation values for air, soft tissue and bone. Considering the average density of each material, at 50 KeV photon energy, the average attenuation values for the air, soft tissue and bone are approximately $2.7\times 10^{-4}$, $0.2$ and $0.7$ cm${}^{-1}$ \cite{hubbell1995tables}. We use two shape dictionaries $\mathfrak{D}^{(1)}$ and $\mathfrak{D}^{(2)}$. The former contains $n_d^{(1)}=5346$ shapes and the latter is in hold of $n_d^{(2)}=6156$ shapes. The basic elements used in the dictionaries are those shown in Fig \ref{fig11}(e). Elements of $\mathfrak{D}^{(1)}$ are chosen to have larger sized shapes as $\phi_1$ is mainly in charge of representing the soft tissue. It is however worth noting that both dictionaries share identical shapes. The shapes are placed all over the imaging domain specially in places that there are chances of objects being present. Plugging the sensitivity relations
\begin{align*}
\frac{\partial \mu}{\partial \alpha^{(1)}_i}=\psi^{(1)}_{\mathcal{S}_i}\Big((\mu_s-\mu_a)\delta_\epsilon(\phi_1)+ (\mu_b-\mu_s)\delta_\epsilon(\phi_1)H_\epsilon (\phi_2) \Big),
\end{align*}
and
\begin{align*}
\frac{\partial \mu}{\partial \alpha^{(2)}_i}=\psi^{(2)}_{\mathcal{S}_i}\Big( (\mu_b-\mu_s)H_\epsilon(\phi_1)\delta_\epsilon (\phi_2) \Big),
\end{align*}
into (\ref{eq56}) and (\ref{eq57}) provides the necessary components in using Algorithm 1 to invert for the level set coefficients. The process of determining a descent direction may be performed in a coordinate descent fashion by alternatively updating the weights for the first and second level set at each iteration. Fig \ref{fig11}(d) shows the result of our inversion. Through a sparse composition of shapes we have been able to reconstruct a reasonable estimate of true image. Considering the image to be a composition of shapes provides a better pose to the problem compared to the case of considering it as an image with piecewise constant regions suitable for a TV inversion.
\subsection{Resistance Tomography: An Archaeological Problem}\label{insec2}
As the last example we consider an archaeological application of imaging subsurface tombs using electrical resistance tomography (ERT). The severely ill-posed nature of the problem as well as the limitations in placing the measuring sensors make this problem very challenging specially when the structures are closely spaced \cite{elwaseif2010quantifying}.

In this technique electric current is injected into the ground and some sensors measure the resulting potential on different regions of the imaging domain. Based on these potential measurements an inverse problem is solved to reconstruct the profile of electrical conductivity. The governing physics may be described as
\begin{align}\nonumber
&\nabla \cdot (\sigma \nabla \rho)=\mathpzc{j} \quad \mbox{in}\; D, \\&\hspace{-.42cm} \xi_1 \sigma \frac{\partial \rho}{\partial \mathbf{n}}+\xi_2 \rho=0\quad \mbox{on} \;\partial D,
\end{align}
where $\mathpzc{j}$ denotes the pattern of injected current, $\sigma$ is the conductivity and $\rho$ is the resulting potential. Functions $\xi_1$ and $\xi_2$ are functions defined on the boundary of the imaging domain and are in charge of imposing appropriate boundary conditions \cite{aghasi2012sensitivity}.

To maintain simplicity, consider $\rho_s$ to be the potential resulted from a point source current $\mathpzc{j}(x)=\delta(x-x_s)$ and $\rho_{s,m}$ representing the measured voltage at points $x_m$ in the domain for $m=1,\cdots M$. The forward model is a nonlinear operator that maps $\sigma(x)$ to the potential measurements. To apply the proposed technique we need to know about the model Jacobian operator. Using the adjoint field technique, it can be shown that perturbations in the measurements are related to the conductivity perturbations through \cite{aghasi2012sensitivity}
\begin{equation}
\delta\rho_{s,m}=\int_D \delta\sigma \nabla \rho_s \cdot \nabla \rho_m \mbox{d}x.
\end{equation}
Here $\rho_m$ is the potential resulted from placing the point source current at $x_m$, known as the adjoint field. When a parametric form $\sigma=\sigma(x,\balpha)$ is considered for the conductivity, we have $\delta \sigma = \sum_i \delta\alpha_i \partial\sigma/\partial \alpha_i$, where $\delta\alpha_i$ is the perturbation of the $i$-th element of $\balpha$. Based on this argument
\begin{equation}
\mathcal{G}_m'(\alpha)\eta= \int_D \nabla \rho_s \cdot \nabla \rho_m (\sum_i \eta_i\frac{\partial \sigma}{\partial \alpha_i}) \mbox{d}x.
\end{equation}
For $\balpha\in\mathbb{R}^{n_d}$, the model Jacobian operator $\mathcal{G}'(\balpha)$ is a matrix of size $M\times n_d$ the $(m,i)$ entry of which is $\int_D \partial\sigma/\partial \alpha_i \nabla \rho_s\cdot\nabla \rho_m\mbox{d}x$. Accordingly, the adjoint operator $\mathcal{G}'^*(\balpha)$ is simply the transposed matrix in this case. This notion can be extended to more complex cases of running different experiments with different sources, as well as more sophisticated source settings such as electric dipoles.

\begin{figure}[!htbp]
\centering
\includegraphics[width=85mm]{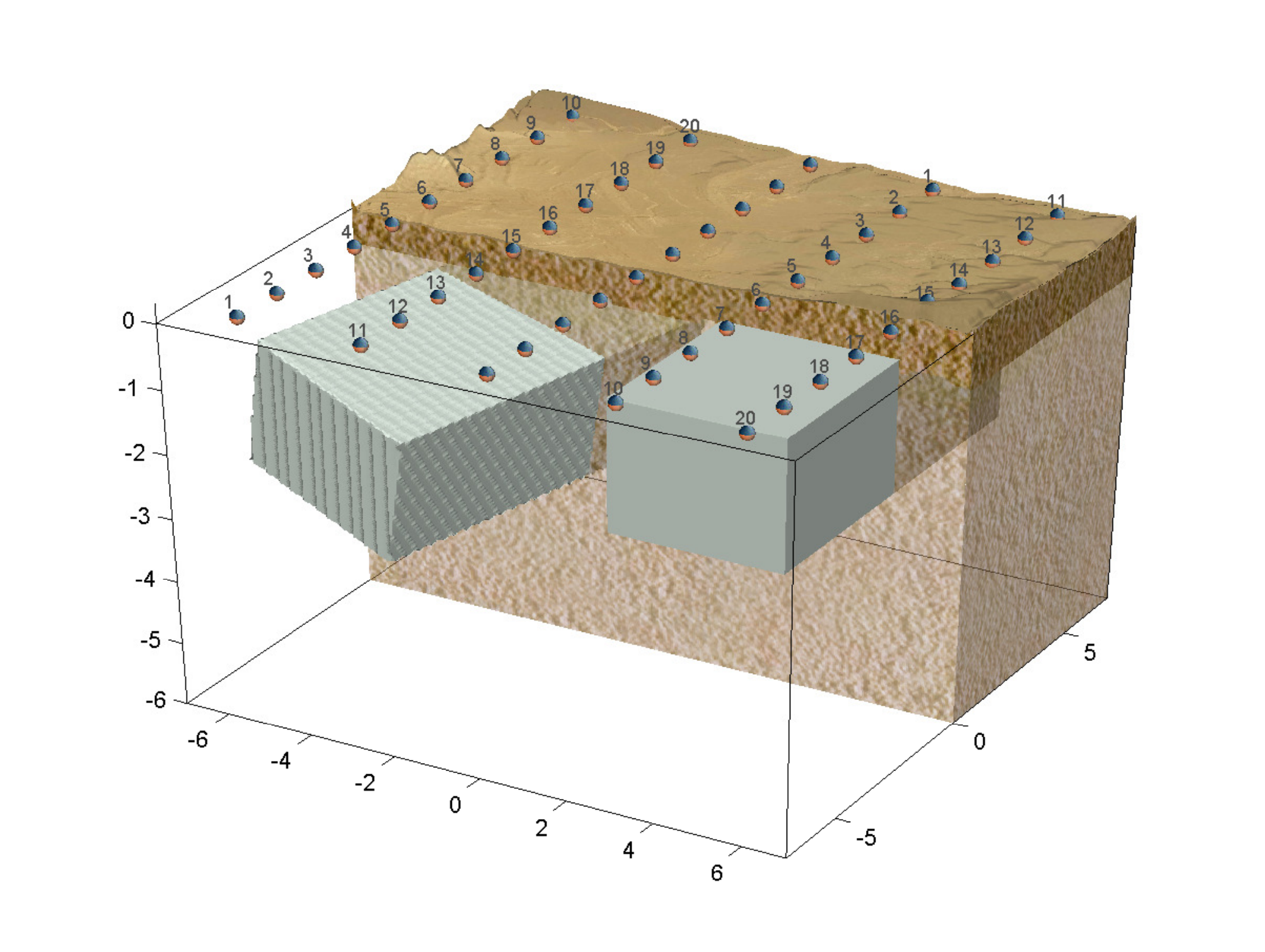}\hspace{-10mm}
\caption{The ERT setting and sensor configuration }\label{fig12}
\end{figure}

\begin{figure}[!htbp]
\centering
\includegraphics[width=140mm]{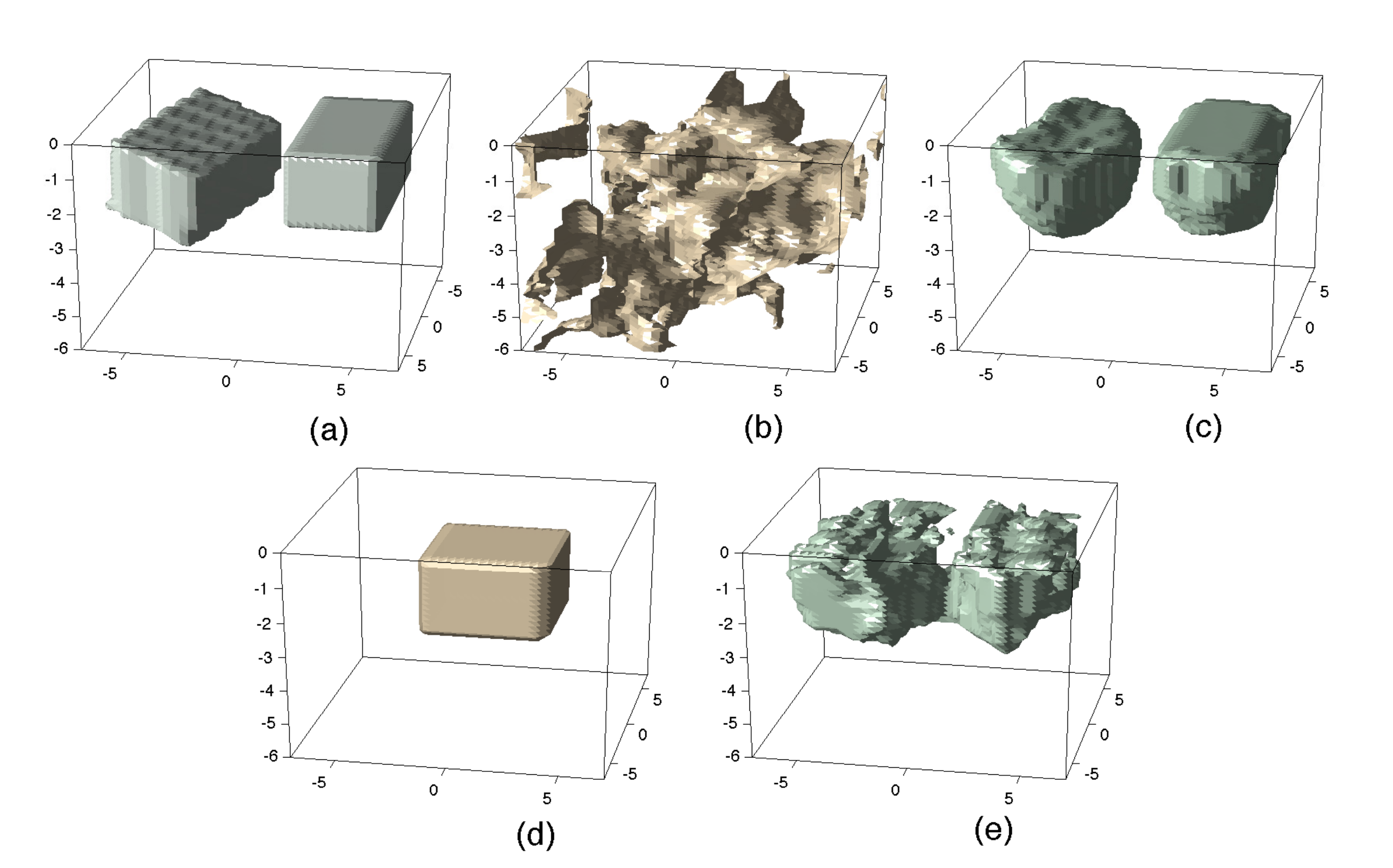}\hspace{-10mm}
\caption{(a) The true tombs (b) Initialization of the algorithm (c) Reconstruction results (d) Initialization of traditional level sets (e) Reconstruction results using traditional level set technique }\label{fig13}
\end{figure}

Fig (\ref{fig12}) shows the setup for a field experiment in a region where two closely placed tombs exist. A total of 50 sensors are placed on the ground and 20 experiments are carried out at each experiment using a pair of sensors as the electric dipole (shown with identical numbers) and the remaining sensors as the measurement sensors. The data is polluted with 1\% white noise. The true cavities are modeled as cubic structures shown in the figure. The structure on the left is slightly tilted in both azimuth and elevation. Following (\ref{eq30}) the conductivity distribution is modeled as
\begin{equation}
\sigma(x,\balpha)= \sigma_s+(\sigma_a-\sigma_s)H_\epsilon\big(\phi(x,\balpha)\big),
\end{equation}
where $\sigma_a$ and $\sigma_s$ are the average known values for the conductivity of air and soil. A total of $n_d=14157$ shape are used in the dictionary. The dictionary elements are cubic structures distributed all over $D$, noting that the true cavities are not among the elements of the dictionary. Fig \ref{fig13}(b) shows the initialization of the algorithm and Fig \ref{fig13}(c) is the imaging results after 11 iterations.

Reconstruction results after 43 iterations using conventional level set technique are also shown in Fig \ref{fig13}(e). We would like to note that the conventional level set inversion is initialized with a very good initialization as shown in Fig \ref{fig13}(e). As observable, the proposed technique is able to provide a more accurate profile of the subsurface conductivity and successfully separate the two structures. The counterpart, however, is unable to provide such level of detail and fails to make the separation in deeper regions of the ground that the sensitivity values are lower.

\section{Conclusion}
The idea presented in this paper may be considered as a new technique in approaching variety of shape-based imaging problems. The main message of this work is changing the geometric problem of shape composition into a variational problem that may be more conveniently analyzed. This conversion is performed rather simply through the notion of pseudo-logical property. In fact following this idea, the shape composition problem becomes similar to the classic problem of representing a function with a weighted sum of functional elements from a dictionary. The notion of sparsity comes next as a means of choosing proper elements, however, the shape-based nature of the problem requires struggling with a nonlinear problem. By considering several examples from different applications, we showed that, although a nonlinear problem, the proposed sparsity promoting technique can successfully handle struggling problems with large dictionary elements. However, more general techniques applicable to larger class of functionals without restrictive assumptions such as knowing the sparsity degree a priori are desirable and welcomed. Of course the non-convex nature of the shape-based problem does not allow talking about uniqueness of reconstructions, however, our group is still working on appropriately restricting the problem for which such analysis is possible. This is still an open arena of research that requires exploring new techniques in nonlinear problems with sparsity constraints.

\appendix
\section{Descent Property of Proposed Steps}
We first show that when $\sigma$ is sufficiently small a direction acquired from (\ref{eq39}) is descent. First a quadratic expansion yields
\begin{equation}\label{ap1}\nonumber
\| \mathcal{G}'(\balpha_k)\bdelta+\mathcal{G}(\balpha_k)\|_\mathbb{S}^2= \|\mathcal{G}(\balpha_k)\|_\mathbb{S}^2+ \je(\balpha)\bdelta + \|\mathcal{G}'(\balpha_k)\bdelta\|_\mathbb{S}^2.
\end{equation}
As $\bdelta_k$ needs to meet the constraint in (\ref{eq39}) we must have $\| \mathcal{G}'(\balpha_k)\bdelta_k+\mathcal{G}(\balpha_k)\|_\mathbb{S}\leq \sigma$, which requires
\begin{align}\label{ap2}\nonumber
\je(\balpha_k)\bdelta_k&\leq \sigma -\|\mathcal{G}(\balpha_k)\|_\mathbb{S}^2- \|\mathcal{G}'(\balpha_k)\bdelta_k\|_\mathbb{S}^2\\&\leq \sigma -\|\mathcal{G}(\balpha_k)\|_\mathbb{S}^2.
\end{align}
The most right side of (\ref{ap2}) is certainly negative for sufficiently small values of $\sigma$, when for instance $\sigma=\min_{\balpha} \|\mathcal{G}(\balpha)\|_\mathbb{S}^2$ or more locally $\sigma< \|\mathcal{G}(\balpha_k)\|_\mathbb{S}^2$. Such choices make $\bdelta_k$ a descent direction for $\mathcal{E}$ at $\balpha_k$.

We next show that when $\tau_{k,\ell}\geq \|\balpha_k\|_1$, direction $\bdelta_{k,\ell}$ acquired from (\ref{eq41}) is descent. Since $ \|\balpha_k+\mathbf{0}\|_1 \leq \tau_{k,\ell}$ and $\bdelta_{k,\ell}$ is a minima for (\ref{eq41}) we need to have
\begin{align}\label{ap3}\nonumber
\| \mathcal{G}'(\balpha_k)\bdelta_{k,\ell}+\mathcal{G}(\balpha_k)\|_\mathbb{S}^2 & \leq \| \mathcal{G}'(\balpha_k)\mathbf{0}+\mathcal{G}(\balpha_k)\|_\mathbb{S}^2 \\&=  \|\mathcal{G}(\balpha_k)\|_\mathbb{S}^2,
\end{align}
which using (\ref{ap1}) simplifies to
\begin{equation}\label{ap4}
\je(\balpha_k)\bdelta_{k,\ell}\leq - \|\mathcal{G}'(\balpha_k)\bdelta_{k,\ell}\|_\mathbb{S}^2\leq 0.
\end{equation}

\section{The Equivalent Problem to the Projection onto an Asymmetric $\ell_1$-Ball}
 \textbf{Proposition:} For $\tilde\balpha\in\mathbb{R}^n$ given, consider $\tilde\balpha_1^\perp$ to be the projection onto an asymmetric $\ell_1$-ball as
 \begin{equation}\label{ap5}
\tilde\balpha_1^\perp= \operatorname*{arg\,min}_{\balpha}\;\;\|\balpha-\tilde\balpha\|_\mathbb{S} \quad s.t. \quad \|\balpha|_{1,w}\leq \tau,
\end{equation}
and $\tilde\balpha_2^\perp$ to be the solution to the weighted $\ell_1$ projection problem
\begin{equation}\label{ap6}
\tilde\balpha_2^\perp=\operatorname*{arg\,min}_{\balpha}\;\;\|\balpha-\tilde\balpha\|_\mathbb{S} \quad s.t. \quad \|\boldsymbol{D}\hspace{-.1mm}\mbox{\small (}\hspace{-.3mm} \small\tilde \balpha \hspace{-.3mm} \mbox{\small )} \normalsize\balpha\|_1\leq \tau,
\end{equation}
where the diagonal matrix $\boldsymbol{D}\hspace{-.1mm}\mbox{\small (}\hspace{-.3mm} \small\tilde \balpha \hspace{-.3mm} \mbox{\small )}$ is defined as (\ref{eq44p5}). Then $\tilde\balpha_1^\perp=\tilde\balpha_2^\perp$.\\

 \textbf{Proof:}
 We first note that both problems (\ref{ap5}) and (\ref{ap6}) are framed as a projection onto a convex set, and in both cases the minima is unique (see Prop B.11 in \cite{bertsekas1999nonlinear}). Clearly elements of $\tilde\balpha_1^\perp$ and $\tilde\balpha$ must have identical signs over the support of $\tilde\balpha_1^\perp$; otherwise by changing the signs of incompatible elements a feasible point is obtained with a lower cost, which contradicts $\tilde\balpha_1^\perp$ being the global minima. A similar argument holds for $\tilde\balpha_2^\perp$ and $\tilde\balpha$, and so
 \begin{equation}\label{ap7}
\boldsymbol{D}\hspace{-.1mm}\mbox{\small (}\hspace{-.3mm} \small \tilde\balpha_1^\perp \hspace{-.3mm} \mbox{\small )} \normalsize\tilde\balpha_1^\perp =\boldsymbol{D}\hspace{-.1mm}\mbox{\small (}\hspace{-.3mm} \small\tilde \balpha \hspace{-.3mm} \mbox{\small )} \normalsize\tilde\balpha_1^\perp,
\end{equation}
and
\begin{equation}\label{ap8}
\boldsymbol{D}\hspace{-.1mm}\mbox{\small (}\hspace{-.3mm} \small \tilde\balpha_2^\perp \hspace{-.3mm} \mbox{\small )} \normalsize\tilde\balpha_2^\perp =\boldsymbol{D}\hspace{-.1mm}\mbox{\small (}\hspace{-.3mm} \small\tilde \balpha \hspace{-.3mm} \mbox{\small )} \normalsize\tilde\balpha_2^\perp.
\end{equation}
Suppose that $\tilde\balpha_1^\perp \neq \tilde\balpha_2^\perp$. Since $\tilde\balpha_2^\perp$ is a solution to (\ref{ap6}), it needs to satisfy $\|\boldsymbol{D}\hspace{-.1mm}\mbox{\small (}\hspace{-.3mm} \small\tilde \tilde\balpha_2^\perp \hspace{-.3mm} \mbox{\small )} \normalsize\tilde\balpha_2^\perp\|_1\leq \tau$. This result along with (\ref{ap8}) and the fact that $\|\balpha|_{1,w}=\|\boldsymbol{D}\hspace{-.1mm}\mbox{\small (}\hspace{-.3mm} \small\balpha \hspace{-.3mm} \mbox{\small )} \normalsize\balpha\|_1$, reveal that $\tilde\balpha_2^\perp$ is also a feasible point for problem (\ref{ap5}). Therefore, based on the strict convexity of the cost
\begin{equation}\label{ap9}
\|\tilde\balpha_1^\perp-\tilde\balpha\|_\mathbb{S}<\|\tilde\balpha_2^\perp-\tilde\balpha\|_\mathbb{S}.
\end{equation}
In a similar fashion equation (\ref{ap7}) reveals that $\tilde\balpha_1^\perp$ is a feasible point for problem (\ref{ap6}) and therefore $\|\tilde\balpha_2^\perp-\tilde\balpha\|_\mathbb{S}<\|\tilde\balpha_1^\perp-\tilde\balpha\|_\mathbb{S}$ which contradicts (\ref{ap9}).

%
%
\end{document}